\documentclass[a4paper,12pt]{elsarticle}
\usepackage{amsmath}
\usepackage[english]{babel}
\usepackage{graphicx}
\usepackage{subfigure}
\usepackage{booktabs}
\usepackage[usenames,dvipsnames]{color} 
\usepackage[numbers]{natbib}
\usepackage[footnotesize,bf,sf,center]{caption}
\usepackage{float}
\usepackage[dvips]{epsfig}
\usepackage[bottom=1.5cm,top=1.5cm,left=1.5cm,right=1.5cm]{geometry}
\usepackage{calc}
\usepackage{gensymb}
\usepackage{soul}
\usepackage{amsmath,amsfonts,amssymb}
\usepackage{setspace}
\usepackage{amssymb}
\usepackage{mathrsfs}
\usepackage{algpseudocode}
\usepackage{dsfont}
\usepackage{epstopdf}
\usepackage{subfig}%
\usepackage[utf8]{inputenc}
\usepackage{amsfonts,amsmath}

\newdefinition{rmk}{Remark}

\usepackage{graphicx}
\newcommand{\vvvert}[1]{{\left\vert\kern-0.25ex\left\vert\kern-0.25ex\left\vert #1 
    \right\vert\kern-0.25ex\right\vert\kern-0.25ex\right\vert}}

\usepackage{multirow}

\newcommand{\cN}{\ensuremath{\mathcal{N}}}

\newcommand{\cS}{\ensuremath{\mathcal{S}}}

\newcommand{\bN}{\ensuremath{\mathbb{N}}}

\newcommand{\bR}{\ensuremath{\mathbb{R}}}

\def\[{\left[}
\def\]{\right]}
\def\<{\langle}
\def\>{\rangle}
\def\({\left(}
\def\){\right)}
\def\[{\left [}
\def\]{\right]}
\def\({\left(}
\def\){\right)}

\newcommand{\norm}[1]{\Vert #1 \Vert}

\DeclareMathOperator*{\argmax}{arg\,max}
\DeclareMathOperator*{\argmin}{arg\,min}

\newcommand{\WW}{\boldsymbol{W}}
\newcommand{\ROM}{\boldsymbol{\Phi}}
\newcommand{\GG}{\boldsymbol{G}}
\newcommand{\SSS}{\Sigma} 
\newcommand{\VV}{\boldsymbol{V}}
\newcommand{\AAA}{\boldsymbol{A}}
\usepackage{tcolorbox}
\usepackage{enumitem}

\newcommand{\vol}{\mathop{\ooalign{\hfil$V$\hfil\cr\kern0.08em--\hfil\cr}}\nolimits}

\usepackage{amsmath} 

\begin{document}
\begin{frontmatter}
\renewcommand\arraystretch{1.0}

\title{\textbf{A fast food-freezing temperature estimation framework using optimally located sensors}}
    
    \author{
   {\bf Felipe Galarce} $^{1}$, \ 
   {\bf Diego R. Rivera} $^{2,3}$, \
   {\bf Douglas R.Q. Pacheco} $^{4}$, \  
   {\bf Alfonso Caiazzo} $^{5}$, \  
   {\bf  Ernesto Castillo$^{2,3}$} \footnote[0]{{\sf Email address:} {\tt ernesto.castillode@usach.cl}, 
     {\sf corresponding author}} \\
     {\small ${}^{1}$ School of Civil Engineering, Pontificia Universidad Católica de Valparaíso, Valparaíso, Chile.} \\
      {\small ${}^{2}$ Departamento de Ingeniería Mecánica, Universidad de Santiago, Santiago, Chile}\\
     {\small ${}^{3}$ Computational Heat and Fluid Flow Lab, Universidad de Santiago de Chile.}\\
     {\small ${}^{4}$ Chair for Computational Analysis of Technical Systems and Chair of Methods for Model-based Development in Computational Engineering, RWTH Aachen University, Germany.}\\
     {\small ${}^{5}$ Weierstrass Institut für Angewandte Analysis  und  Stochastik,  Leibniz-Institut  im  Forschungsverbund  Berlin  e.V. (WIAS)}
    }


\begin{abstract}
This article presents and assesses a framework for estimating temperature fields in real time for food-freezing applications, significantly reducing computational load while ensuring accurate temperature monitoring, which represents a promising technological tool for optimizing and controlling food engineering processes. The strategy is based on (i) a mathematical model of a convection-dominated problem coupling thermal convection and turbulence, and (ii) a least-squares approach for solving the inverse data assimilation problem, regularized by projecting the governing dynamics onto a reduced-order model (ROM). The unsteady freezing process considers a salmon slice in a freezer cabinet, modeled with temperature-dependent thermophysical properties. The forward problem is approximated using a third-order WENO finite volume solver, including an optimized second-order backward scheme for time discretization. We employ our data assimilation framework to reconstruct the temperature field based on a limited number of sensors and to estimate temperature distributions within frozen food. 
Sensor placement is optimized using a novel greedy algorithm, which maximizes the observability of the reduced-order dynamics for a fixed set of sensors. The proposed approach allows efficient extrapolation from external sensor measurements to the internal temperature of the food under realistic turbulent flow conditions, which is crucial for maintaining food quality.
\end{abstract}
\end{frontmatter}


\section{Introduction}

Freezing processes are crucial for food transport and preservation \cite{sridhar2021food,hu2022novel} and have important economic implications, relying on large amounts of energy to satisfy safe freezing conditions. 
Food freezing involves a highly nonlinear liquid-solid phase change that produces significant variation in processed food and strongly depends on the operational conditions imposed, for example, in a freezing tunnel. 
Liquid-solid phase change implies a boundary evolution scenario in which the phases interact, natural convection occurs, the thermophysical properties of the food change strongly, and turbulence develops in the airflow that freezes the food \cite{GONZALEZ2021101101,DASILVA202246,du2024physics}. In addition, food humidity, texture, color and flavor depend on the evolution of the freezing rate \cite{CFDfood2023,szpicer2023application,ouyang2025improvement}. All these factors influence food quality, and the possibility of estimating in advance how the food
temperature changes during freezing is of both scientific and technological interest.

Computational fluid dynamics has become a standard tool to support the understanding of the physics of fluids and of the heat-mass transfer that occurs in food processing. For example, computer simulations have been successfully used to evaluate convective coefficients \cite{HU2001702,MORAGA202176,ray2024framework}, to estimate freezing times \cite{HOANG202093,piechnik2024experimentally,xu2025effect}, and to test the effectivity of different operational conditions \cite{pr12091852}, thus improving overall thermal processes. Numerical methods used for these simulations have been largely based on finite difference \cite{WANG2007502,FERREIRA2020326}, finite element \cite{sujinda2025simulation,khan2024computational,zaman2024thermal} or finite volume methods \cite{MORAGA200317,MULOT201913}. 
The strengths and weaknesses of each of these numerical methodologies have been recently discussed in \cite{takhar2017-review1,muktiarni2023-review2,arpaci2024exploring}. 

Inverse problems in the food industry have gained attention in recent years, primarily focusing on parameter estimation techniques. Recent studies, such as those in \cite{REDDY2022110909}, emphasize
the use of inverse strategies, particularly in identifying material properties and diffusivity parameters in food products. For example, research in \cite{FABBRI201463} highlights the estimation
of moisture diffusivity by inverting a finite element model applied to various solid foods through the Levenberg–-Marquardt algorithm. 
Just as ensuring the properties of heat-treated foods is of scientific and technological interest,
efficient energy consumption in the processes is also an aspect that can be optimized using similar techniques \cite{ravula2023estimation}. In  \cite{CORNEJO201637}, comparison with experimental data showed that the bootstrap technique can estimate both specific heat and thermal conductivity in a freezing range.

In this work, we address the estimation of the temperature fields in a freezing salmon slice from measured data, focusing on the case of \textit{partial} and \textit{low-resolution} observations. This type of research belongs to the broader field of \textit{data assimilation}, where available data and mathematical models are used both to estimate unknown parameters of a PDE
(\textit{parameter estimation}) \cite{YADAV2015182} or to fully discover model fields (\textit{state estimation}) \cite{DUDA2016201, shaoqing_2020,PHAM2024109801}. A variety of data-driven methods have emerged recently in the literature, powered up by the gathering of large data-sets. Nevertheless, this kind of approach, though undoubtedly capable, may suffer drawbacks in scenarios where there is much less access to repetitive measurements, or whenever physical interpretability is a must \cite{alzubaidi2023survey}. This calls for alternatives that find a compromise between the purely data-driven approach and the physically informed approach, thus exploiting some well-known governing laws of the underlying phenomena. \cite{HAIK2023115868, wang2024_causality, Toscano2025, Brunton2025, Garay2024}. In the latter case, as in the problem considered in this work, the goal is to estimate time- and space-resolved fields over the whole computational domain or on a subset of interest, starting from experimental measurements that are typically noisy. This usually leads to ill-conditioned problems, and the literature propose several strategies to convey such estimations, such as maximally weighted approaches \cite{YU2023108169} or classical Kalman filters \cite{SHEN2024122452}.

For the forward modeling, we use a finite volume solver with a third-order weighted essentially non-oscillatory method (WENO3) for the cell-face fluxes and an implicit optimized second-order backward differentiation formulation (BDF2-opt) to handle the temporal evolution, combined with the URANS $k$-$\omega$ SST turbulence model. The efficient implementation of these methods enables accurate and fast approximation of turbulent thermally-coupled flows and the nonlinear phase change process. However, the small space and time scales required to describe the phase-change phenomenon \cite{HAJJAR2020105243}, together with the need to systematically repeat these simulations---as needed for the inverse problem---calls for using Reduced-Order Modeling (ROM) strategies. With that, we can reliably decrease the dimensionality of the problem while still preserving its fundamental characteristics. Well-established ROM techniques include the reduced-basis method \cite{virginie_ReducedBasis}, proper orthogonal decompositions \cite{SVDantoulas}, randomized decompositions for large scale systems \cite{rSVD2024}, modern hybrid approaches based on machine learning techniques \cite{ROMOR2025113915},  high order dynamic mode decompositions \cite{CORROCHANO2023108219}, non-linear approaches \cite{MAZZILLI2022106915}, among many others (see e.g. \cite{frangos2010, benner2017model} for a review on ROM strategies).

As an inverse technique, we employ a regularized least-squares approach based on the so-called Parametrized Background Data Weak (PBDW) method, see, e.g.,  \cite{HAIK2023115868, maday19_adaptPBDW, MPPY2015,cohen_2020}. 
In this approach, the optimization problem to be solved
for the state estimation is based on the projection of the
physics on a reduced-dimensional space that must reproduce the observed data in a weak form. Application of the PBDW in biomedical problems has been recently presented in \cite{GGLM2021,GLM2021,galarce2023_MRE,GMC2024}. In this work, we use a methodology inspired by the PBDW to assimilate the temperature fields in the liquid-solid phase change, including turbulent flow conditions. In particular, we employ the ROM-based component of PBDW, which is equivalent to the full PBDW approach with infinite model regularization, as explained in \cite{HAIK2023115868}, assuming no bias on the physical model. Additionally, we combine the state estimation with a Greedy algorithm for the optimal placement of temperature sensors. Greedy algorithms have been used for experimental design in the past years. In \cite{greedy1}, a Bayesian design is used for seismic full waveform inversion. In \cite{greedy2}, a greedy algorithm was employed using a minimal variance error for parameter estimation. In \cite{greedy3}, optimal sensors are computed targeting directly a mean squared error between the model and measurements. Greedy methods have been combined with the PBDW in, e.g., \cite{manoharSensorPlacement, binevGreedy,HAIK2023115868}. A recent related work is \cite{mula2025}, where the problem is posed so that an optimal stabilization of the inverse approach is acquired and applied in a Hamiltonian system for wave-type phenomena. In \cite{mula2025_was}, sparse measurement locations are sought from a dictionary of available positions using a Wasserstein barycenter criteria. In contrast with these works, our proposal solely relies on available \textit{a priori} results for the regularized least-squares reconstruction, prioritizing sensor locations that maximize the \textit{observability} of a static reduced-order model. 

In summary, this article proposes a temperature estimation framework from optimally placed synthetic data and a regularized least-squares estimator through a physics-informed projection, using turbulent simulations of the airflow surrounding the freezing food. The temperature estimation algorithm is proven valid for extrapolating fields in non-observable regions, keeping good reconstruction errors due to a properly designed sensor placement strategy. As a further application of the fully reconstructed field, we also show the capabilities of the algorithm to predict temperature-related quantities, such as the local freezing curve and the local freezing rate, aspects of high interest in food engineering. This proof-of-concept article will specifically use synthetically generated measurements from a thermal camera. 

The rest of this manuscript is organized as follows. Section \ref{Sec_ProblemS} states the problem, including the physical situation studied in this work, the governing equations solved and the models used, and the initial and boundary conditions.  Section \ref{Sec_DirectA} summarizes the finite volume method used as a direct strategy, focusing on the distinctive aspects with respect to standard finite volume solvers used in the literature for food process modeling, and a standard validation of the method with well-established numerical and experimental benchmarks. In Section \ref{Sec_Inverse}, an inverse pipeline for state estimation is presented, including the reduced-order strategy used and a fully devoted sub-section where an optimal sensor placement strategy is adapted for the problem at hand. In Section \ref{Sec_NumEX}, the feasibility of the proposed inverse method to estimate states in processes involving food freezing is evaluated exhaustively using different sensor positions and highlighting the representativeness of the methodology for quantities of interest to the food industry.  We finally draw concluding remarks in Section \ref{Sec_Conclusion}.

\section{Mathematical modeling of the conjugate food freezing process} \label{Sec_ProblemS}

In this work, we consider a two-dimensional (2D) setup for the forward simulations and the inverse reconstruction, as depicted in Fig. \ref{fig:physicalsituation}$(a)$. That geometrical simplification is based on the symmetry of the three-dimensional problem. We also performed a validation with experimental data previously published, which confirms that the use of the 2D model is adequate for the physical situation studied \cite{RIVERA2024110558}. In that article, the importance of using ad-hoc turbulence models to describe the freezing evolution was discussed in detail, with URANS strategies being sufficiently accurate and efficient to reproduce experimental results. Since the key novelty of this article is the inverse estimation, we accept the 2D direct solution and use the 2D results as the ground truth values. The inverse methodology, without the inclusion of the greedy algorithm, was systematically tested on 3D scenarios with good outcomes in \cite{galarce2023_MRE,GGLM2021}. In this concern, we do not restrict the analysis to the two dimensional case due to methodological limitations, since there are none, but rather because we would like to introduce a complete minimal reconstruction pipeline for a realistic scenario, still with a 2D set-up.

\subsection{Model setting} \label{modelsetings}
Even though humidity and ice formation are important for freezing applications \cite{jia2022control,sun2023regulating,ribeiro2024effects}, including the mass transfer equation and more sophisticated phase-change models to characterize the airflow-food interaction is highly challenging from the forward modeling viewpoint. Also, including thermodynamic relations to compute the water condensation when the vapor concentration is higher than the saturated vapor concentration is necessary to upgrade the physical modeling of freezing applications. Since this work is focused on inverse modeling and temperature state estimation, with the finite volume solutions as input data, the mathematical model defined to represent the salmon freezing process only considers the momentum, continuity, and energy conservation laws. The particular treatment used to define the airflow and the food is defined separately in the following subsections. 

\begin{figure}[h]
    \centering
    \includegraphics[width=1\linewidth]{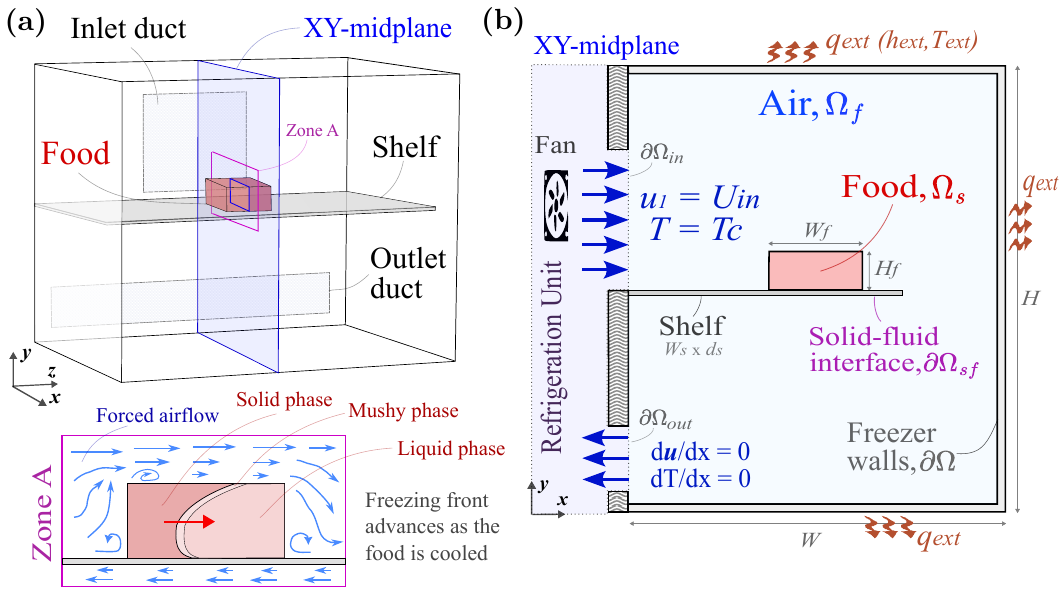}
    \caption{(a) Sketch of the physical domain for the food freezing process with a brief description of the freezing phenomena. (b) Details of the two-dimensional section considered in this study.}
    \label{fig:physicalsituation}
\end{figure}

We model the freezing process of food inside a convective freezer that circulates cold air between the cabinet and the refrigeration unit through inlet and outlet ducts. The conjugate dynamics of the process results from the interplay of forced heat convection in the turbulent airflow and the nonlinear heat transfer coupled with the liquid-to-solid phase change process within the food. This setup is illustrated in Fig.~\ref{fig:physicalsituation}. 

The freezer is assumed to have a cabinet of width $W=0.35$ m, height $H=0.4$ m, and depth $D=0.45$ m. Cold air coming from the refrigeration unit is forced into the cabinet by a fan at a constant velocity $U_{in}$ and temperature $T_{C}$ through an inlet duct of height $H_{in}=0.12$ m, and it leaves the cabinet through an outlet duct of height $H_{out}=0.06$ m. An idealized hexahedral food slab of width $W_{f}=0.08$ m, height $H_{f}=0.03$ m, and depth $D_{f}=0.08$ m is placed on a glass shelf of width $W_{s}=0.25$ m and thickness $d_{s}=0.002$ m, located at the mid-height of the cabinet. Also, we assume that the walls of the cabinet exchange heat with the environment by conduction through the insulated layers and by natural convection with the surrounding air. The rear wall that separates the cabinet from the refrigeration unit is considered adiabatic. In what follows, we will denote with $\Omega_f \subset \mathbb R^2$ and $\Omega_s \subset \mathbb R^2$ the fluid-flow and food domains, respectively, and $[0,t_f]$ is the physical time interval considered for the simulation.

\subsection{Airflow dynamics} \label{airflowdynamics}

The airflow in the fluid domain is assumed to be described by the incompressible unsteady Reynolds-averaged Navier--Stokes equations (URANS) for velocity 
$\boldsymbol{u}$
and pressure $p$, coupled to a heat convection equation for the air temperature $T$. Using a Boussinesq approximation for the buoyancy forces to deal with the thermal coupling, the airflow obeys the following governing equations: 
\begin{equation}
\begin{aligned}
    \nabla \cdot \boldsymbol{u} = 0 \qquad \text{in} \ \Omega_f \times (0,t_{f}\rbrack, \\
        \rho \partial_{t} \boldsymbol{u} + \rho\boldsymbol{u} \cdot \nabla \boldsymbol{u} - \nabla \cdot [2(\mu + \mu_{t}) \nabla^s \boldsymbol{u}] + \nabla p = \boldsymbol{f} \qquad \text{in} \ \Omega_f \times (0,t_{f}\rbrack, \\
        \rho \partial_{t} T + \rho\boldsymbol{u} \cdot \nabla T  - \nabla \cdot[(\lambda /C_{p} + \sigma_{T}\mu_{t}) \nabla T]  = 0 \qquad \text{in} \ \Omega_f \times (0,t_{f}\rbrack.   
\end{aligned}
\label{eq:gov}
\end{equation}

In Eq.~\eqref{eq:gov} $\rho$, $\lambda$, $C_{p}$, and $\mu$ denote the density, thermal conductivity, specific heat, and viscosity of the air, respectively. These physical properties are assumed as constant at a reference temperature of $T_{ref}=0$ $\degree$C. This reference temperature was determined as the average between the maximum and minimum temperature values during the process. Then, Ideal Gas and Sutherland’s Laws were used to compute the values of density (1.2930 kg/m$^{3}$), thermal conductivity (0.02436 W/mK$^{1}$), and viscosity (1.3316 $\times$ $10^{-5}$ Ns/m$^{2}$), while the specific heat was assumed as constant with a value of 1003.68 J/(kgK) and the thermal expansion coefficient was computed as $T_{ref}^{-1}$.

The source term in the momentum equation is defined as:
$$
\boldsymbol{f} = \rho\beta_{ref} \left( T - T_{ref} \right)\boldsymbol{g} + \frac23 \nabla \cdot (\rho k),
$$
and includes the buoyancy and the turbulent source terms, where $\beta_{ref}$ is the air's thermal expansion coefficient and $k$ is the turbulent kinetic energy. The additional turbulent viscosity $\mu_{t}$ introduced by the URANS method is described by the $k$-$\omega$ SST model. 
The model attempts to predict turbulence and close the problem using two partial differential equations for two variables: $k$ and $\omega$. The first variable is the turbulence kinetic energy, while the second is the specific rate of dissipation. The turbulence model is based on two transport equations defined as:
\begin{equation} \label{eq:turb}
\begin{aligned}
        \rho \partial_{t} k + \rho\boldsymbol{u} \cdot \nabla k - \nabla \cdot (\mu + \sigma_{k}\mu_{t}) \nabla k = \tilde{P}_{k}-\beta^{*}\rho k \omega \qquad \text{in}  \ \Omega_{f} \times(0,t_{f}\rbrack,
        \\
        \rho \partial_{t} \omega + \rho\boldsymbol{u} \cdot \nabla \omega  - \nabla \cdot (\mu + \sigma_{\omega}\mu_{t}) \nabla \omega  = \frac{\rho \gamma}{\mu_{t}}\tilde{P}_{k}-\rho \beta \omega^{2} + C_{k\omega} \qquad \text{in}  \ \Omega_{f} \times (0,t_{f}\rbrack.  
\end{aligned}
\end{equation}

By using the two turbulence equations, the turbulent viscosity is defined 
as a function of the turbulent kinetic energy and specific turbulence dissipation rate ($\omega$) as:
\begin{equation}\label{eq:mu_t}
\mu_{t} = \rho \,\frac{a^{*} k}{\max(a_{1}\omega , SF_{2})} \,.
\end{equation}

In Eq.~\eqref{eq:mu_t}, $S=\sqrt{2\nabla^s \boldsymbol{u}:\nabla^s \boldsymbol{u}}$, $a^{*} = 0.31$, and $F_2$ is an adjustable function of the model.
In Eq.~\eqref{eq:turb}, $\beta^{*}=0.09$, $\tilde{P}_{k}= \min (\tau : \nabla\boldsymbol{u} ,10\beta^{*}\rho\omega k)$, $C_{k \omega}= 2\sigma_{\omega,2} \omega^{-1} (1-F_{1}) (\nabla k \cdot \nabla\omega)$
and $\sigma_{T}$, $\sigma_{k}$, $\sigma_{\omega}$, $\beta$, and $\gamma$ are functions varying as a convex combination of inner and outer constant values (close to the high gradient regions and far from them) with the blending function $F_{1}$, e.g.,  $\sigma_{\cdot}=F_{1} \sigma_{\cdot,1} +(1-F_{1}) \sigma_{\cdot,2}$, where the indices $1$ and $2$ denote the inner and outer values, respectively. 
For the detailed model, we refer the reader to \cite{MenterkwSST2003}.

\subsection{Food freezing process} \label{foodfreezingprocess}

The internal convection of the liquid water contained in the muscle fibers and connective tissue of salmon can be neglected at the macro-level \cite{moraga1999numerical,tabilo2024freezing, isafahan2024}. To account for the effect of latent heat during the phase change process in the salmon slab, the thermal properties (the density $\rho_s$, the specific heat capacity $C_s$, and the thermal conductivity $\lambda_s$) are considered as temperature-dependent, according to the effective heat capacity method \cite{comini1974,CUI2022101593}.

The heat transfer and the liquid-solid phase change process in the food slab are described by a non-linear energy equation, assuming an incompressible and homogeneous material neglecting mass transfer, such that:
\begin{equation}
    \begin{aligned}
    \partial_{t} (\rho_{s} C_{s} T_{s})  - \nabla \cdot (\lambda_{s} \nabla T_{s})  = 0 \qquad \text{in}  \ \Omega_{s} \times (0,t_{f}\rbrack.    
    \end{aligned}
    \label{eq:govfood}
\end{equation}

\begin{figure}
    \centering
    \includegraphics[width=1\linewidth]{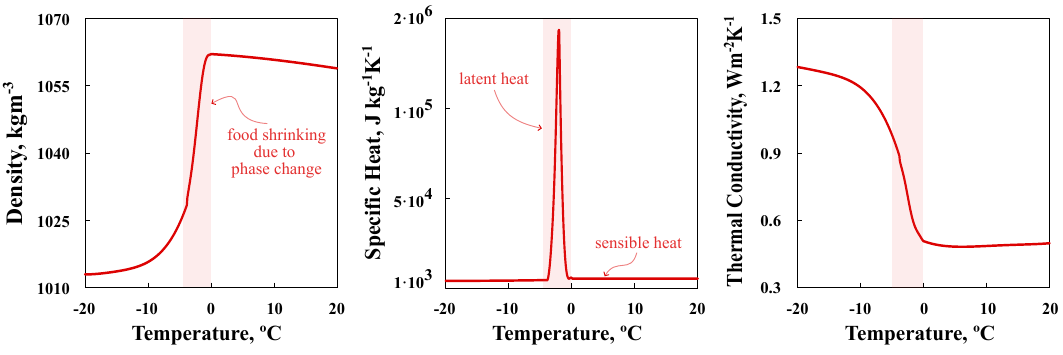}
    \caption{Temperature dependency of the thermal properties of salmon meat. Left: density ($\rho_s$), center: specific heat capacity ($C_s$), right: thermal conductivity ($\lambda_s$).}
    \label{fig:foodprops}
\end{figure}

The temperature variability of the properties is depicted in Fig.~\ref{fig:foodprops}, assuming third-order fitting with (piecewise-constant) temperature-dependent coefficients:
\begin{equation}
\psi = \sum_{j=0}^3 a^{\psi}_j(T) T^j,    
\end{equation}
for $\psi = (\rho_{s}C_{s}$,$\lambda_s$).
The values of the coefficients for the volumetric heat capacity $\rho_s C_s$ and thermal conductivity $\lambda_s$ are presented in Table \ref{tab:salmonprops}. 

\begin{table}[h]
    \centering
    \renewcommand{\arraystretch}{1.5} 
    \begin{tabular}{c c c c c c}
        \hline
        Property & $T\in$ & [-25°C, -5°C] & (-5°C, -3.5°C] & (-3.5°C, 0°C] & (0°C, 25°C]  \\
        \hline
        \multirow{4}{4em}{$\rho_{s}C_{s}$} & $a_{0}$ & 4.298$\times10^{7}$ & 3.046$\times10^{10}$ & -4.604$\times10^{6}$ & 4.365$\times10^{6}$ \\
        & $a_{1}$ & 7.166$\times10^{6}$ & 1.906$\times10^{10}$ & -1.037$\times10^{8}$ & -3.228$\times10^{5}$ \\
        & $a_{2}$ & 4.199$\times10^{5}$ & 4.002$\times10^{9}$ & -1.095$\times10^{8}$ & 2.535$\times10^{4}$ \\
        & $a_{3}$ & 8.031$\times10^{3}$ & 2.817$\times10^{8}$ & -4.260$\times10^{7}$ & -5.838$\times10^{2}$ \\
        \hline
        \multirow{4}{4em}{$\lambda_{s}$} & $a_{0}$ & 5.654$\times10^{-1}$ & -3.276$\times10^{-1}$ & 5.202$\times10^{-1}$ & 5.122$\times10^{-1}$ \\
        & $a_{1}$ & -9.014$\times10^{-2}$ & -5.611$\times10^{-1}$ & 1.181$\times10^{-2}$ & -3.900$\times10^{-3}$ \\
        & $a_{2}$ & -4.038$\times10^{-3}$ & -9.160$\times10^{-2}$ & 2.707$\times10^{-2}$ & 3.585$\times10^{-4}$ \\
        & $a_{3}$ & -6.721$\times10^{-5}$ & -5.767$\times10^{-3}$ & 1.081$\times10^{-3}$ & -7.767$\times10^{-6}$ \\
        \hline
    \end{tabular}
    \caption{Coefficients for the polynomial functions of volumetric heat capacity ($\rho_{s}C_{s}$) and thermal conductivity ($\lambda_{s}$) of salmon for the different temperature ranges.}
    \label{tab:salmonprops}
\end{table}

\subsection{Boundary, initial, and interface conditions} \label{boundaryconditions}

Boundary conditions at the inlet $\partial\Omega_{in}$ are prescribed using constant air velocity $U_{in}$ and constant temperature $T_{C}$, which are assumptions based on: (i) the airflow enters horizontally with an average velocity that correlates with typical mass flow rates of medium-sized fans, and (ii) the refrigeration unit rapidly cools the air to the set freezing point during most of the charging process. 

At the outlet $\partial\Omega_{out}$, we assume a developed flow condition, where velocity is multiplied by a correction factor to ensure mass conservation: 
$F_{m}=\int\boldsymbol{u}\ d\Omega_{in}\ /\int\boldsymbol{u}\ d\Omega_{out}$. At the solid-fluid interfaces $\partial\Omega_{sf}$, non-slip boundary conditions are imposed for velocity, heat flux continuity for temperature, and model-defined values are specified for the turbulent variables. For the internal walls of the freezer $\partial\Omega$, we computed a Robin-type thermal boundary condition with an equivalent thermal resistance: $R_{W}=1/h_{ext} + \delta x_{w}/\lambda_{w}$. $R_{W}$ includes the heat conduction resistance by the insulated layers of the freezer walls, with $\lambda_{w}=0.026$  Wm$^{-1}$K$^{-1}$ and $\delta x_{w}=0.05$ m \cite{MARQUES2013457}, and the heat convection resistance with the external air at temperature $T_{ext}$, assuming typical correlations of the external convective coefficient $h_{ext}$ with Nusselt and Rayleigh numbers in the form of $h_{ext}=\lambda Nu(Ra)/L_{c}$ \cite{ghajar2014heat}. The boundary conditions are summarized in Table \ref{tab:boundaryc}.
\begin{table}[h]
    \centering
    \renewcommand{\arraystretch}{1.5} 
        \begin{tabular}{c c c c c} 
        \hline
         $\Omega$ & Velocity & Temperature & Kinetic e. & Dissip.~rate of $k$ \\
        \hline
         $\partial\Omega_{in}$: & $\boldsymbol{u}=(U_{in},0)$ & $T=T_{C}$ & $k=(0.02 U_{in})^{2}$ & $\omega = \sqrt{k}/\beta^{*}H_{in}$ \\
        \hline
         $\partial\Omega_{out}$: & $ \nabla\boldsymbol{u}\cdot \boldsymbol{n}=0$ & $\nabla T\cdot \boldsymbol{n}=0$ & $\nabla k\cdot \boldsymbol{n}=0$ & $\nabla \omega\cdot \boldsymbol{n}=0$ \\ 
        \hline
        $\partial\Omega$: & $\boldsymbol{u}=0$ & $\lambda \nabla T \cdot \boldsymbol{n}=(T_{ext}-T_{w})/R_{W} $ & $k=0$ & $\omega = 800\mu/\rho d_{wall}^{2}$ \\
        \hline
        $\partial\Omega_{sf}$: & $\boldsymbol{u}=0$ & $\lambda \nabla T \cdot \boldsymbol{n}= \lambda_{s}\nabla T_{s} \cdot \boldsymbol{n}$ & $k=0$ & $\omega = 800\mu/\rho d_{wall}^{2}$ \\ 
        \hline
        \end{tabular}
    \caption{Boundary conditions at the inlet duct $\partial\Omega_{in}$, the outlet duct $\partial\Omega_{out}$, the freezer walls $\partial\Omega$, and at the solid-fluid interface $\partial\Omega_{sf}$; $d_{wall}$ is the distance from the walls to the closest cell centroid.}
    \label{tab:boundaryc}
\end{table} 

Initial conditions in the fluid domain $\Omega_{f}$ are given by zero flow velocity, zero pressure field, and zero values of the turbulent variables, i.e. $\boldsymbol{u}_{0}=0$, $p_0=0$, $k_{0}=0$, and $\omega_{0}=0$. Also, the fluid and solid domains are in thermal equilibrium with the environment, i.e., $T_{0}={T_s}_0=T_{ext}$.

The inverse pipeline of the present work requires a dataset of snapshots from a set of time-dependent direct simulations with a variable space of parameters. Therefore, the magnitudes of the main boundary conditions are varied within a fixed range: $U_{in} \in [0.15,0.25]$ \text{m/s}, $T_{C} \in [-18,-26] \ ^{\circ}\text{C}$, $T_{ext} \in [18,26] \ ^{\circ}\text{C}$, and $h_{ext} \in [0.4,1.2] \ \text{W/}\text{m}^{2}\text{K}$. 
The detailed procedure for constructing this five-dimensional set is explained in Section \ref{Sec_NumEX}.

Since natural convection is important in the problem at hand, we need to define the Rayleigh number for each studied case: 
\begin{align*}
  Ra = \frac{|\boldsymbol{g}| \beta (T_{C} - T_{ext})L_{c}^{3}}{\nu \alpha} ,  
\end{align*}
where for each case $T_C$ and $T_{ext}$ change and the thermophysical properties must be redefined using $T_C$. The characteristic length is defined as the height of the freezer cabinet in all cases. Using this Rayleigh definition, we can confirm that the values are very high ($Ra>10^{10}$), which is the reason for including a turbulence model in the direct simulations.

\section{Direct problem}\label{Sec_DirectA}

This section is devoted to describing the numerical approach we adopt to approximate solutions of the differential equations \eqref{eq:govfood}. We select a finite volume scheme, introduced next in detail, along with the temporal discretization of the problem. Next, we test our forward solver with a beef freezing set-up, benchmarking the numerical outcomes with available localized experimental measurements.

\subsection{Finite volume method} \label{sec_discretization}

To compute a numerical solution, the system of equations \eqref{eq:gov}, \eqref{eq:turb},\eqref{eq:govfood} is discretized in space using the finite volume method (FVM) by decomposing the computational domain into $N_{v}$ non-overlapping cells $\Omega_{i}$. Integrating the PDEs over each cell and applying the divergence theorem yields:
\begin{equation}
    \int_{\Omega_{i}} \left( \rho \partial_{t}\phi \right) \mathrm{d}\Omega_{i} = \int_{S_i} \left(\rho\boldsymbol{u}\phi \right) \cdot \boldsymbol{n}\, \mathrm{d}S_i + \int_{S_i} \left(\Gamma \nabla\phi \right) \cdot \boldsymbol{n} \, \mathrm{d}S_i + \int_{\Omega_{i}} \boldsymbol{f} \, \mathrm{d}\Omega_{i},
    \label{eq:fvm}
\end{equation}
where $\phi$ stands for any state variable ($\boldsymbol{u},T,k,\omega$), $\Gamma$ is the corresponding diffusion coefficient, $\boldsymbol{f}$ stands for the corresponding forcing term, $S_i$ is the cell boundary, and $\boldsymbol{n}$ is the outward normal unit vector to $S_i$. We use a structured anisotropic Cartesian mesh with a staggered grid approach that stores the cell-averaged value of scalars at the cell centroid and velocities at the cell faces. 

We selected the Weighted Essentially Non-Oscillatory (WENO) scheme because of its advantages over other commonly used higher-order schemes such as the Total Variation Diminishing (TVD) family, from both numerical and  implementation standpoints. The WENO scheme achieves arbitrarily high-order accuracy in smooth regions by combining nonlinear weights and lower-order stencils, maintaining stable and non-oscillatory solutions. Furthermore, using a monotonic flux function (Godunov or Lax-Friedrichs, for example), spurious oscillations (undershoots/overshoots) are eliminated, ensuring monotonicity-preserving methods \cite{barrenechea2024finite,bozorgpour2025recent}. This issue is particularly important when sharp gradients in temperature or velocity occur. Also, implementing the WENO method is straightforward in anisotropic Cartesian structured meshes, a key property for mechanical science applications.
The cell face reconstruction used is based on a third-order WENO method (WENO3) \cite{Shu1998} that combines two adaptive weights, $w_{0},w_{1} \ge 0$, with $w_{0}+w_{1}=1$, and two-point sub-stencils, $S_{i}^0(\phi_{i-1},\phi_{i})$ and $S_i^{1}(\phi_{i},\phi_{i+1})$, such that: $\phi_{S_i} = w_{0}S_i^{0} + w_{1}S_i^{1}$. For incompressible turbulent flows such as those studied here, including dissipation is a must, and the key feature is the correct dissipative structure of the method and its ability to represent physically admissible solutions, which WENO3 fulfills \cite{tsoutsanis2014weno,SHU2016598,navas2024exploring,hao2024large,jayswal2025smooth}, as discussed in the direct model validation section. For the WENO3 approach, each sub-stencil approximates $\phi_{S_i}$ using a linear interpolation through the respective stencil points (cell centroids) in the anisotropic mesh. 
The weights are computed as a function of the state variable as:
\begin{equation*}
w_k = \frac{\alpha_k}{\alpha_0 + \alpha_1}\ ;\ \alpha_k = \frac{d_k}{\left(\beta_k + \epsilon\right)^2}\ ;\ \beta_{k}=\left(\phi_{i-1+k}-\phi_{i+k}\right)^2, \qquad \text{for}\ k=0 \ \text{and} \ 1,
\end{equation*}
where $d_{0}=\frac13$, $d_{1}=\frac23$, and $\epsilon=10^{-6}$. Finally, we use the Lax--Friedrichs procedure to ensure a consistent monotone physical flux at the cell faces, and we use a second-order central difference method for the gradient reconstruction at cell faces. 

The time discretization is based on a uniform time step dividing the considered time interval 
$\left\lbrack 0,t_{f} \right\rbrack$ into $N_{t}$ sub-intervals of size $\mathrm{\Delta}t = t^{n + 1} - t^{n}$. The state variables on the right-hand side of \eqref{eq:fvm} are treated implicitly, and the time derivative is approximated using an optimized backward differentiation formula (BDF2-opt). This method is unconditionally stable and uses four previous time-levels, such that:
\begin{equation*}
    \partial_{t}\phi = \Delta t^{-1}\left(c_{1}\phi^{n+1} + c_{2}\phi^{n} + c_{3}\phi^{n-1} + c_{4}\phi^{n-2}\right),
\end{equation*}
where the coefficients are calculated by a linear combination of a classical second-order and a third-order BDF schemes \cite{Vatsa2010,ORTEGA2021345}, such that: $ c_i = \chi c_i^{\textrm{BDF2}} + (1-\chi)c_i^{\textrm{BDF3}}$. In all simulations considered in this, we set $\chi=0.52$. 

Velocity and pressure are coupled using the segregated pressure-correction method SIMPLEC, and the set of resulting linearized algebraic equations is solved by the Bi-Conjugated Stabilized Gradient method (Bi-CGSTAB) with the ILU preconditioner. 
The SIMPLEC algorithm offers a significant advantage over other segregation methods, such as SIMPLE or PISO, particularly when addressing highly nonlinear problems. This advantage stems from its lower computational cost per outer iteration, as only one inner iteration is needed, which enhances the overall nonlinear convergence and the coupling between the unknown fields. The problem being solved involves the coupling of turbulent viscosity and buoyancy forces in the linear momentum equation, as well as temperature-dependent properties of food in the energy equation. This study effectively used a fixed-point iterative strategy to linearize the problem, making it straightforward to implement within the SIMPLEC segregation method.
The maximum relative difference of all state variables between two successive iterations is used as converge criteria. The solution converges when that error is less than $10^{-6}$.

Regarding the boundary conditions, Dirichlet-type ones are imposed on the finite volume nodes, and Neumann boundary conditions are defined using first-order finite differences considering non-uniform grids. Robin conditions are also computed with first-order finite differences. Near-wall grid refinement using a tangent hyperbolic function was essential to deal with the imposition of turbulent boundary conditions. Grid refinement allows us to define the first cell height (which dimensionless is $y^+$) close to the recommended value for URANS turbulent models, $y^+=\mu_{t} \delta x_n/\nu \approx 1$.

\subsection{Verification of the finite volume solver} \label{sec_validationfvm}

To verify the numerical method, we solve two well-established benchmark problems: 
\begin{itemize}
    \item[(P1)] Freezing process of beef inside a conventional freezer cabinet by turbulent natural convection of air, with a Rayleigh number of $Ra=6.81\times 10^{7}$ \cite{MORAGA200317}, and
    \item[(P2)] Turbulent natural convection of air in a square cavity differentially heated on the sides, with  $Ra=1.58\times 10^{9}$ \cite{AMPOFO20033551}.
\end{itemize}

\begin{figure}
    \centering
    \includegraphics[width=1\linewidth]{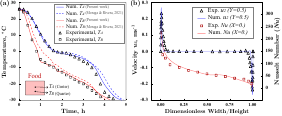}
    \caption{Numerical results for the present direct solver (continuous line), for a previous direct solver \cite{MORAGA202176} (dashed line), and experimental data (squares and triangles) for the two benchmark problems: (a) P1: Temperature evolution at two selected points inside the food, as illustrated in lower-left corner. (b) P2: Profiles of vertical velocity and Nusselt number at the mid-height and at the hot wall of the cavity, respectively.}
    \label{fig:validationdp}
\end{figure}

Fig.~\ref{fig:validationdp}$(a)$ compares the experimental data with present and previous numerical results from the authors (from \cite{MORAGA202176}) for P1. The collection of experimental data was obtained using two well-calibrated 0.3 mm diameter T-type thermocouples, which recorded the temperature every 10 seconds at two fixed points within the food (at center and quarter height, see Fig.~\ref{fig:validationdp}$(a)$) with a standard uncertainty of 0.5$\ ^{\circ}\text{C}$. The present work uses an improved numerical solver compared to \cite{MORAGA202176}, with high-order schemes (WENO3 and BDF2-opt) and the $k$-$\omega$ SST turbulence model that improves the accuracy of results. Considering these methods and the careful discretization in space and time of the present numerical procedure, we estimate that the numerical/experimental uncertainty ratio ($UR$) is less than 1, $0.2\le UR\le1$.

For P2, Fig.~\ref{fig:validationdp}$(b)$ compares the numerical results for the considered configuration with the experimental profiles of vertical air velocity at the mid-height of the cavity $(Y=0.5)$ and the Nusselt number at the hot wall $(X=0)$. Despite the thin boundary layer, the results agree with the experimental data, obtaining RMS error values of 0.05 ms$^{-1}$ and 18.3 for the velocity and Nusselt numbers, respectively. From these and the previous results, we conclude that the considered Finite Volume Method can accurately describe the turbulent convection of air and the food freezing process.

\subsection{Forward numerical simulations} \label{Sec_Direct}

In this Section we discuss the main features of the forward numerical results related to the forced heat convection within a turbulent airflow and to the phase change process of the food. We use a representative parameter set defined as the mean values of the range for each boundary condition, i.e.,  $U_{in}=0.2$ ms$^{-1}$, $T_{C}=-22 \ ^{\circ}\text{C}$, $T_{ext} =22 \ ^{\circ}\text{C}$, and $h_{ext}=0.8  \ \text{W}\text{m}^{-2}\text{K}^{-1}$. 

To define the number of volume cells and the time step size required to ensure the independence of solutions regarding discretization, we conduct a convergence study by comparing velocity and temperature profiles along control lines. We define three different meshes composed of $18,410$ (M1 mesh), $22,802$ (M2 mesh), and $29,070$ (M3 mesh) cells. In this work, we use a wall-concentrated anisotropic Cartesian mesh defined to capture velocity and temperature gradients and accurately capture turbulent quantities. The meshes used were constructed by dividing the domain into eighteen subdomains to control the size of the control volumes separately. The minimum and maximum wall-normal lengths dimensions used to construct each mesh are defined as $h_{min,i}$ and $h_{max,i}$, being $i$, a counter that refers to the specific mesh. In specific, $h_{min,1}= 4.0 \cdot 10^{-3}$ and $h_{max,1}= 8.5 \cdot 10^{-2}$, $h_{min,2}=7.2 \cdot 10^{-3}$ and $h_{max,2}= 1.8 \cdot 10^{-3}$, and $h_{min,3}=1.3 \cdot 10^{-4}$ and $h_{max,3}= 5.1 \cdot 10^{-3}$. Concerning the discretization, we compare three different time step sizes of $\Delta t_1=0.03$ s, $\Delta t_2=0.05$ s, and $\Delta t_3=0.08$ s, and we evaluate the effect of this discretization parameter in the numerical results. 

\begin{figure}[h]
    \centering
    \includegraphics[width=1\linewidth]{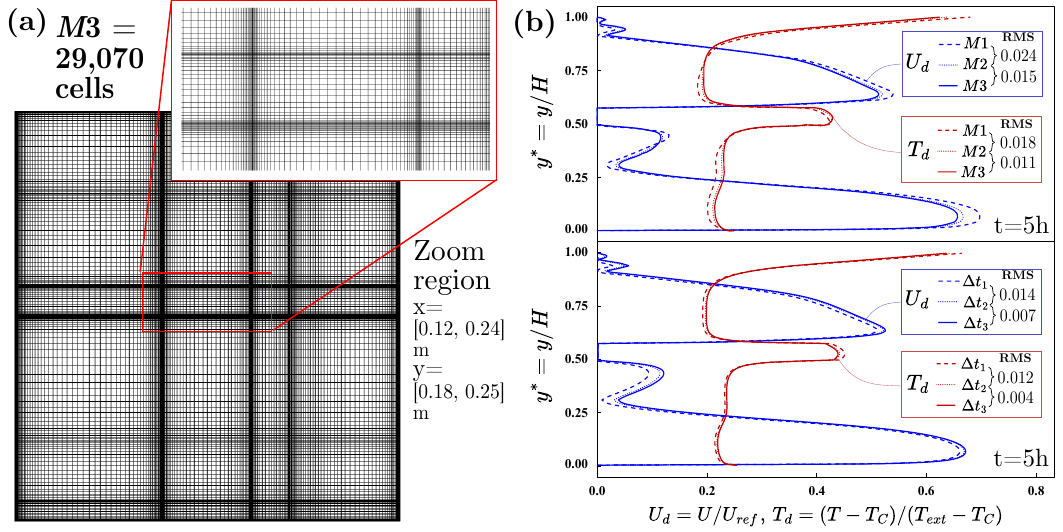}
    \caption{(a) Mesh refinement of the selected mesh $M3$ with a zoom region near the food. (b) Discretization study: comparison of velocity and temperature profiles between three mesh sizes (upper) and three time steps (lower).}
    \label{fig:discretization_study}
\end{figure}

In Fig. \ref{fig:discretization_study}(b), we compare the dimensionless vertical profiles of velocity and temperature at the freezer cabinet mid-plane at $t=5$ h (hours), using $U_{ref}=0.5$ ms$^{-1}$. When comparing the root mean square error (RMS) at fixed points between each successive mesh size and time step, we found that the results changed on average by 2.1\% from $M1$ to $M2$, and by 1.3\% from $M2$ to $M3$. Regarding the time step, the results varied by 1.4\% from $\Delta t_1$ to $\Delta t_2$, and by 0.5\% from $\Delta t_2$ to $\Delta t_3$. Given the low variation of the results with the discretization parameters, we selected the mesh $M3$ and the time step $\Delta t_2$ to conduct the numerical test through the work.

\begin{figure}[h]
    \centering
    \includegraphics[width=1\linewidth]{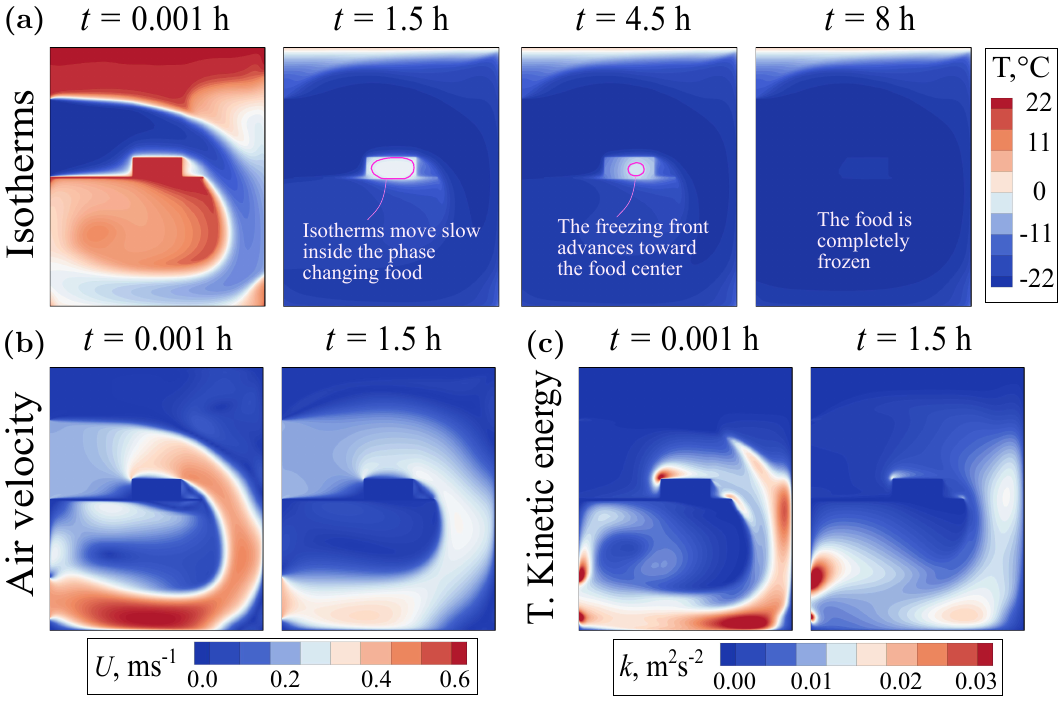}
    \caption{ Depictions of different snapshots of the temperature fields $(T)$, air velocity fields $(U)$, and turbulent kinetic energy fields $(k)$ for the representative set of parameters.}
    \label{fig:fig1direct}
\end{figure}

Fig.~\ref{fig:fig1direct} illustrates snapshots of the temperature, air velocity as $U=\sqrt{u_{1}^2+u_{2}^2}$, and turbulent kinetic energy fields at characteristic time instants: $t=0.001$ h, $t=1.5$ h, $t=4.5$ h, and $t=8$ h. The results indicate that the inlet airflow causes rapid cooling of the air domain within the first minute, with some high-temperature stratifying isotherms in the upper region produced by natural convection effects. Contrary, the food cools down slowly due to its high thermal inertial caused by the latent heat during phase change. Moreover, both air velocity and turbulent kinetic energy fields indicate that at the very beginning the flow is more unsteady and chaotic, evidencing that the fluctuations and gradients are notably higher. However, the flow stabilizes rapidly within the first hour and reaches a state of periodic oscillations of velocity and temperature that are unnoticeable since they occur on a small time scale.

\begin{figure}
    \centering
    \includegraphics[width=1\linewidth]{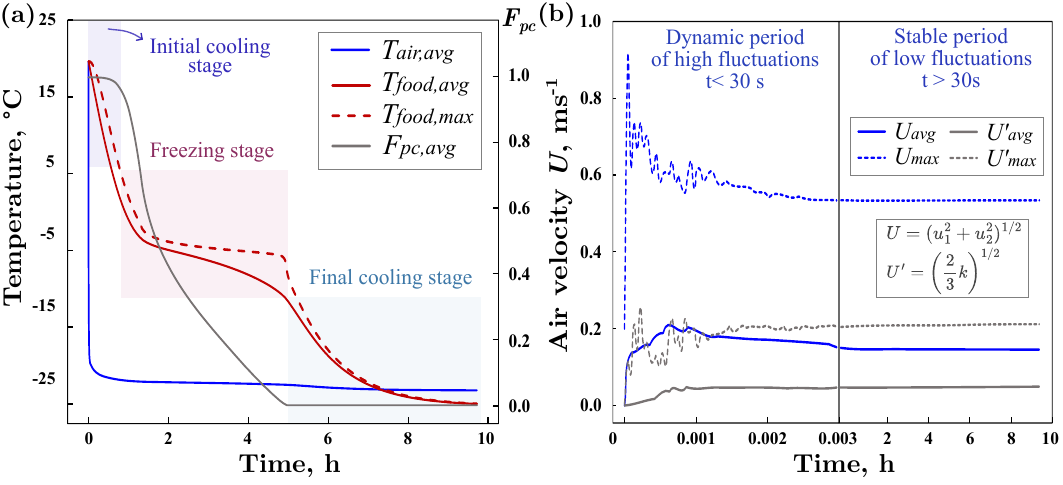}
    \caption{Left: Evolution of volume-averaged and maximum values of temperature in the air and food, and of the food liquid fraction. Right: evolution of volume-averaged and maximum values of the time-averaged air velocity $U$ and of the instantaneous fluctuation of velocity $U'$.}
    \label{fig:fig2direct}
\end{figure}

To understand the different time scales of the fluid dynamics and heat transfer processes, Fig.~\ref{fig:fig2direct} depicts the time evolution of volume-average and maximum values of temperatures (air and food), of food liquid fraction ($F_{pc}$), and of time-averaged and instantaneous fluctuations of velocity (i.e. $U$ and $U'=\sqrt{(2/3)k}$), revealing the abrupt drop in air temperatures. In contrast, the average and maximum temperature in the food illustrates a characteristic freezing curve where the effect of the latent heat on the temperature-slope is clearly defined in the freezing stage. 
Regarding the evolution of volume-average and maximum values the two components of the air velocity in a URANS model, the most dynamic patterns of $U$ and $U'$ are observed at the first seconds of the freezing process, but then these values quickly stabilize until they reach a steady behavior over time. Moreover, the elevated value of $U_{max}$ is about four times higher than its volume-average value, demonstrating the presence of high-velocity gradients in the flow domain. Lastly, the relation between the instantaneous fluctuation and the time-average velocity, or turbulent intensity $I_{t}=U'/U$, reaches volume-average and maximum values of 22\% and 38\%, respectively, indicating the convective-dominant nature of the problem. These points evidence the need for a suitable turbulent model that accurately computes the increased momentum and heat diffusion due to the turbulent mixing phenomenon.

\begin{figure}[h]
    \centering
    \includegraphics[width=1\linewidth]{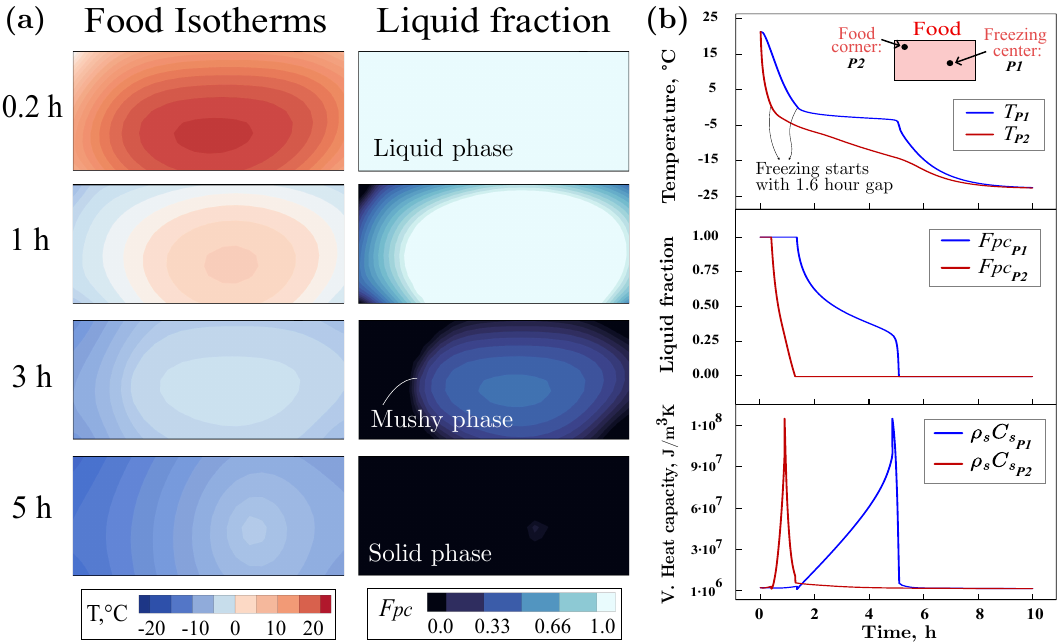}
    \caption{Left: Snapshots of temperature and liquid fraction fields within the food domain. Right: local evolution of the temperature, liquid fraction, and volumetric heat capacity at two control points.}
    \label{fig:fig3direct}
\end{figure}

Fig.~\ref{fig:fig3direct} (left) shows snapshots of the temperature and liquid fraction fields within the food at $t=0.2$ h, $t=1$ h, $t=3$ h, and $t=5$ h. As expected, the liquid fraction field evolution is directly correlated with the isotherms since it is a temperature-derived quantity. Another relevant aspect observed in the results is that the freezing center moves during the freezing process. Initially, the high-temperature isotherms are concentrated just below the geometric center, and during the phase change, they are rather concentrated close to the right surface of the food. This behavior occurs because the cold airflow slips along the upper wall and hits the food on the left, while the lower right corner is located precisely in a flow recirculation region, which weakens heat transfer. 

To highlight the spatial differences in the dynamics of the food freezing process, Fig.~\ref{fig:fig3direct} (right) depicts the temporal evolution of temperature, liquid fraction, and volumetric heat capacity at two fixed points within the food domain. These points are conveniently selected as the freezing center ($P1$), and in the upper corner region ($P2$) where the freezing rate is expected to be highest. The difference in the path of temperature and liquid fraction evolution between the two points within the food shows the relevance of the locality in the phase change phenomenon. The interior of the food experiences a more isothermal process than the periphery, indicating a lower heat transfer rate. Consequently, the maximum values of volumetric heat capacity related to the latent heat occur more slowly and at different periods of time in the freezing process.

\section{Inverse problem}\label{Sec_Inverse}

We aim to reconstruct the temperature field from the available data by solving an optimization problem that minimizes the discrepancy between the model and the measurements. Specifically, the problem is formulated as follows: Given an experimental dataset $\ell \in \bR^m$, in a fixed time, where $m$ is the number of measurements (pixels in a picture, for instance), find a field $T^*$ (e.g., the reconstructed temperature) such that:
\begin{equation}
T^* = \argmin_{T \in \bR^\cN} \norm{\ell - \WW^T T}_2.
\label{eq:ls}
\end{equation}
The minimizer $T^*$ is sought in an $\cN$-dimensional space, where $\cN$ is the number of degrees of freedom from the numerical solution considered for the forward problem (proportional to the number of cells $N_v$ in the finite volume discretization), and where $\WW \in \bR^{\cN \times m}$ is a suitable matrix representing the observation process, herein assumed linear on the state variable $T$.

The main proposal of this work is a procedure for estimating temperature fields from experimental data. In this context, there are several alternatives for handling the input data, which we can summarize in two types of measurements: Point-wise measurements and averaged temperature sensors in low-resolution temperature fields from thermo-cameras. In any case, the sensor gives local averages of the temperature field to be estimated, and therefore, the linearity assumption on the state to be estimated is reasonable.

As schematically depicted in the right bottom corner of Fig.~\ref{fig:pipeline}, measuring a high dimensional temperature state is understood as applying the observer matrix $\WW$ on a temperature state. The matrix columns, otherwise understood as the Riesz representers of the measurements \cite{galarce2023_MRE}, are  indicatrix functions on the data pixels, i.e., for every thermo-camera pixel $\Omega_i \subset \Omega$ we have:
$$
\ell_i(\Omega_i; x) = \left\lbrace \begin{aligned}
    1 / \norm{\ell_i(\Omega_i)}, &\quad x \in \Omega_i, \\
    0, &\quad x \in \Omega \backslash \Omega_i,
\end{aligned}
\right.
$$
for $i=1,\ldots,m$ pixel measurements.



\subsection{Least-squares reconstruction with ROM regularization}

Typically, the number of available measurements is smaller than the dimension of the numerical solution ($\cN > m$). As a consequence, problem \eqref{eq:ls} does not have, in general, a unique solution. This can be seen by developing the optimal conditions of the optimization problem and noticing that the resulting linear system, i.e., 
\begin{equation}
\WW \WW^T T = \WW \ell,
\end{equation}
is not invertible, as the projection matrix $\WW \WW^T$ has at most $m$ non-zero eigenvalues. Otherwise stated, for the usual case $\cN > m$, there are infinitely many temperature fields whose measurements are equal. To overcome this issue, the strategy adopted in this work relies on seeking the minimizer in a suitable subspace of  dimension $n \ll \cN$ (ROM - reduced-order model) defined by the basis:
\begin{equation*}
\ROM = \left( \phi_1 | \ldots | \phi_n \right) \in \bR^{\cN \times n}\,.
\end{equation*}
Hence, the minimization problem \eqref{eq:ls} is replaced by:
\begin{equation}
c^* = \argmin_{c \in \bR^n} \norm{\ell - \WW^T \ROM c}_2.
\label{eq:ls_rom}    
\end{equation}

The optimality conditions for the reduced-order problem \eqref{eq:ls_rom} yield the well-known \textsl{normal equations} \cite{galarce2023_MRE}:
\begin{equation}
    \GG^T \GG c = \GG^T \ell,
    \label{eq:normal_eqs}
\end{equation}
where $\GG:=\WW^T \ROM$ is a cross-Gramian matrix. Then, the temperature field in the original space can be recovered doing $T^* = \ROM c^*$.

Properly built, the reduced-order space defined by $\ROM$ will maintain the most relevant features of the original transport problem described by Eq.~\eqref{eq:gov}. A practical approach to defining a suitable physics-based basis $\ROM$ will be described in more detail in Section \ref{ssec:rom}.

The dimensional reduction has two important implications:
\begin{itemize}
    \item Speed-up: The reduced model in $\ROM$  has a considerably smaller dimension than the original, full-order problem. Provided that this low-dimensional space can represent sufficiently well the physics of the original problem, the overall complexity is considerably reduced, since only an $n \times n$ system shall be inverted. 
    \item Physics-driven recovery: The reduced-order model is tightly coupled to the original physics. Hence, it regularizes the otherwise non-invertible problem \eqref{eq:ls} by reducing the number of unknowns with a physics-informed search space for the optimization problem.
\end{itemize}

The analogous of problem \eqref{eq:ls} in a continuous and infinite-dimensional setting has been studied, from a general point of view, in \cite{MPPY2015, BCDDPW2017, cohen2022_rom}. The regularized least-squares approach used in this article is one of the numerical ingredients of the Parametrized Background Data-Weak (PBDW) method proposed in \cite{MPPY2015}, or it can be seen as the perfect-model (without bias) case explained in \cite{HAIK2023115868}. The analysis of the resulting low-dimensional problem reveals that the condition $\cN > m > n$—indicating that the number of observations exceeds the dimension of the reduced-order model—is necessary to ensure well-posedness. This condition will be assumed valid throughout the remainder of this article.


Thus, the reconstruction problem is considered to have two stages. An offline (or training) stage dedicated to the simulation of several direct problems, plus the computation of the ROM and the observation matrix $\WW$. In the offline stage, most of the computational cost is focused, in contrast with the online phase, where the sole inversion of a $n\times n$ system of equations enables real-time computations.

\subsection{Reduced-order modeling}
\label{ssec:rom}
This section discusses the generation of the reduced-order model $\ROM$ using a parametrized singular value decomposition (SVD) of a set of snapshots from the forward problem
(named in certain contexts PCA -- principal component analysis -- or POD -- proper orthogonal decomposition -- in its  generalization to arbitrary infinite-dimensional functional settings). The goal is to define a low-dimensional space that can sufficiently well approximate the physics of the phase change problem and the turbulent heat convection developed by the airflow. 
For this purpose, the first step consists of concatenating snapshots 
\begin{equation*}
\AAA = \left(
T(\theta_1) | \hdots 
 | T(\theta_k) \right),
\end{equation*}
of the numerical solution of Eq.~\eqref{eq:gov} for a (large) amount of configurations of the parameters defining the underlying model $\theta_1,\hdots,\theta_K$. In the above notation, the vector $\theta_i \in \bR^p$ ($p \in \bN$) defines a generic set of parameter values for \eqref{eq:gov}, e.g., time, fan configurations, boundary conditions, phase-change parameters, food-specific coefficients, among others.

The basis of the reduced-order space $\ROM$ is then computed using the unique and truncated singular value decomposition of the snapshot matrix $\AAA$, namely:
\begin{equation}
    \label{eq:POD}
    \AAA = \ROM\SSS \VV^T =
    \underbrace{
    \begin{bmatrix}
    \textcolor{red}{\blacksquare} & \textcolor{red}{\blacksquare} & \cdots & \textcolor{red}{\blacksquare} \\
    \textcolor{red}{\blacksquare} & \textcolor{red}{\blacksquare} & \cdots & \textcolor{red}{\blacksquare} \\
    \vdots & \vdots & \ddots & \vdots \\
    \textcolor{red}{\blacksquare} & \textcolor{red}{\blacksquare} & \cdots & \textcolor{red}{\blacksquare} \\
    \textcolor{red}{\blacksquare} & \textcolor{red}{\blacksquare} & \cdots & \textcolor{red}{\blacksquare}
    \end{bmatrix}
    }_{\ROM \ (\mathcal{N} \times n)}
    \underbrace{
    \begin{bmatrix}
    \sigma_1 & 0 & \cdots & 0 \\
    0 & \sigma_2 & \cdots & 0 \\
    \vdots & \vdots & \ddots & \sigma_n
    \end{bmatrix}
    }_{\Sigma \ (n \times n)}
    \underbrace{
    \begin{bmatrix}
    \textcolor{green}{\blacksquare} & \textcolor{green}{\blacksquare} & \cdots & \textcolor{green}{\blacksquare} & \textcolor{green}{\blacksquare} \\
    \textcolor{green}{\blacksquare} & \textcolor{green}{\blacksquare} & \cdots & \textcolor{green}{\blacksquare} & \textcolor{green}{\blacksquare} \\
    \vdots & \vdots & \ddots & \vdots & \vdots
    \end{bmatrix}^T
    }_{\VV^T \ (n \times k)},
\end{equation}
defining $\ROM$ using the first $n$ columns of the full decomposition. We will call $\sigma_i$ the diagonal entry $\SSS_{ii}$ of the matrix $\SSS$. We see that the core idea of the SVD regularization for the inverse technique relies upon constraining our optimization search to a linear combination of the columns in $\ROM$ (red squares in \eqref{eq:POD}).

\subsection{Greedy algorithm for optimal sensor placement}
\label{sec:greedy}

As shown in, e.g., \cite{GGLM2021,MPPY2015}, the following \textsl{a-priori} error bound for the reconstructed temperature field $T^*$ holds 
\begin{equation}\label{eq:Tstart-estim}
\frac{\norm{T - T^*}}{T} \leq \left( \inf_{c \in \bR^n} \frac{\norm{\GG c}}{\norm{\ROM c}} \right)^{-1} \frac{\norm{A - A_n}_{\text{FRO}}}{\norm{A}_{\text{FRO}}},
\end{equation}
where $\AAA_n$ is the $n$-th low rank approximation of $\AAA$ and ${\norm{\cdot}_{\text{FRO}}}$ denotes the Frobenius norm:
\begin{align*}
    {\norm{A}^2_{\text{FRO}}} = \mathrm{tr}(A^TA) = \sum_{i,j}A_{ij}^2\, .
\end{align*}
The bound \eqref{eq:Tstart-estim} shows where the condition $m > n$ comes from, indicating that a reduced model that is too enriched would lead to non-observable modes, a critical issue for well-posedness in data assimilation problems. This provides a criterion for selecting the modes to be used in the inverse method, guarantying that any other configuration for the matrix $\ROM$, would lead to sub-optimal reconstruction results (this does not mean that $\ROM$ always approximates better the forward dynamics as $n$ increases).
Using standard arguments, the second term in the right-hand-side of \eqref{eq:Tstart-estim} can be expressed in terms of the singular values of $\AAA$:
$$
\frac{\norm{A - A_n}_{\text{FRO}}}{\norm{A}_{\text{FRO}}} = \frac{\sqrt{\sum_{i=n+1}^{\text{max}\{ \cN, K \}} \sigma_{i}^2}}  {\sqrt{\sum_{i=1}^{\text{max}\{ \cN, K \}} \sigma_{i}^2}}.
$$
For the first term, let us consider the spectral induced $\ell^2$ norm, considering the singular value decomposition $\GG = \hat{U} \hat{S} \hat{V}^T$, and as usual, denote with $\hat{S}_i$ the entry $\hat{S}_{ii}$ of the matrix $\hat{S}$. 
By construction, the columns of $\hat{V}$ are a basis of $\bR^n$, allowing us to write, for any $c \in \mathbb R^n$, $c = V d$ for a unique $d \in \bR^n$. 
Since both the columns of $V$ and of $\ROM$ are orthonormal, we obtain
\begin{equation}\label{eq:GPhi}
\begin{aligned}
\left(\frac{\norm{\GG c}}{\norm{\ROM c}} \right)^2 & =
\frac{ d^T \hat{V}^T \GG^T \GG \hat{V} d} {d^T V^T \ROM^T \ROM V d}
= \frac{ d^T \hat{S}^2 d} {d^T d}
=\frac{ \sum_{i=1}^n d_i d_j \hat{S}_{i}^2} {\sum_{i=1}^n d_i^2}, \\
& \ge \frac{ \hat{S}_{n}^2 \sum_{i=1}^n d_i^2 } {\sum_{i=1}^n d_i^2} = \hat{S}_{n}^2,\;
\forall c \in \mathbb R^n,
\end{aligned}
\end{equation}
where we have also used the fact that the singular values are decreasing and positive. 
The inequality \eqref{eq:GPhi} allows to conclude that
$$
\inf_{c \in \bR^n}
\frac{\norm{\GG c}}{\norm{\ROM c}} \ge \hat{S}_{n},
$$
and using \eqref{eq:Tstart-estim} to obtain the following a priori bound for a given dimension $n$ of the reduced-order space
\begin{equation}
e(n) = \hat{S}_{n}^{-1} \frac{\sqrt{\sum_{i=n+1}^{\text{max}\{ \cN, K \}} \sigma_{i}^2}}  {\sqrt{\sum_{i=1}^{\text{max}\{ \cN, K \}} \sigma_{i}^2}}.
\label{eq:apriori}
\end{equation}

Equation \eqref{eq:apriori} has two main implications for the properties of the temperature reconstruction. On the one hand, better approximation properties
of the reduced-order model lower the enumerator of the error bound.
On the other hand, the term $\hat{S}_n^{-1}$ (the last singular value of the cross-Gramian matrix $\GG$) indicates that 
increasing the dimension of the reduced-order space, i.e., decreasing $\hat{S}_n^{-1}$,
lowers the quality of the expected result.
This trade-off can be interpreted as a lack of
\textit{observability} of the reduced-order space if the dimension is too large.
From the point of view of linear algebra, the observability can be defined as the angle between the subspace spanned by the columns of the matrix $\WW$ and the reduced-order model spanned by the columns of $\ROM$ \cite{BCDDPW2017}. 

We will use \eqref{eq:apriori} as a criterion to improve the ROM observability and to compute the optimal placement of measurement sensors.
In the rest of this section, we describe an iterative approach to systematically increase the dimension of the matrix $\WW$ (i.e., the number of measurements) by locating the sensors in such a way as to maximize the observability (i.e., minimizing \eqref{eq:apriori}) for a fixed dimension $n$ of the reduced-order model.

Using a criteria based on \eqref{eq:apriori} to select the optimal dimension of the ROM and data assimilation has been already applied in other contexts \cite{GGLM2021,GLM2021,galarce2023_MRE}. In \cite{manoharSensorPlacement}, a classical QR pivoting strategy is employed (thus limiting the algorithm to square systems of equations), whereas in \cite{HAIK2023115868}, a criteria based on \eqref{eq:apriori} is used, nonetheless with certain differences: First, with our method, we do select the optimal measurements based on a fixed dimension for the ROM (which could be necessary for certain engineering applications, where a certain level of approximation should be ensured for the reduced basis). Second, in \cite{HAIK2023115868}, the measurement addition is based upon optimizing the \textsl{least-stable} mode iteratively instead of using all modes at once, as we propose here. Both approaches may deliver similar results as $m$ increases, but for a fixed amount of measurements (as in many engineering scenarios), we compute the global stability using the full reduced basis, thus ensuring that the conditioning is optimal for fixed $n$ and $m$.

\begin{tcolorbox}[colback=blue!5!white, colframe=blue!75!black, title=\textbf{Algorithm 1: Greedy Sensor Placement}, fonttitle=\bfseries, sharp corners]
\begin{itemize}[leftmargin=4em]
    \item[\textbf{Input:}] 
        \begin{itemize}
            \item ROM basis $\ROM$ (and its fixed dimension $n$).
            \item Pool of feasible sensor positions $\cS$.
        \end{itemize}
    \item [\textbf{Step 1}] Select the first sensor as:
    $\hat{\ell}_1 := \argmax_{w \in \mathcal{S}} \norm{w^T \ROM}.$
    
    \item [\textbf{Step 2}] Select the next $n-1$ sensors (necessary to ensure well-posedness) as:
    $$ 
    \hat{\ell}_i := \argmax_{w \in \mathcal{S}}  \sigma_i (\hat{\mathbf{W}}(w)^T \ROM) ,\; \text{for } i=2,\dots,n, 
    $$
    where 
    $ \hat{\mathbf{W}}(w) := \big[ \hat{\ell}_1 \mid \cdots \mid \hat{\ell}_{i-1} \mid w \big] $
    is the $\mathcal{N} \times i$ observation matrix for a candidate $i$-th sensor $w$, and $\sigma_i(B)$ denotes the $i$-th singular value of any matrix $B$.

    \item [\textbf{Step 3}] For $i > n$, locate any additional sensor maximizing the smallest non-zero singular value of the resulting sensor matrix, i.e.,
    $$ 
    \hat{\ell}_i := \argmax_{w \in \mathcal{S}} \sigma_n (\hat{\mathbf{W}}(w)^T \ROM) . 
    $$
\end{itemize}
\end{tcolorbox}

Then, for a given $n$, we aim at selecting the location of $m \geq n$ sensors from a pool of $M$ possible choices over the observable region $\Omega$. The proposed approach follows
the underlying idea of a \textit{Greedy} method, adding the sensors progressively and looking at each stage for 
the location of a new sensor that yields the smallest a-priori bound \eqref{eq:apriori}. Let us recall that, as seen in \ref{Sec_Inverse}, a
measurement sensor can be seen as a linear operation on the temperature field. Thus, considering the Riesz representation of this linear functional, in what follows, we will refer to a \textit{sensor} as a column of the matrix $\WW$. Moreover, we will denote as 
$\cS = \{ \ell_1, \ldots, \ell_M \}$ a set of all possible sensor  locations, and $\hat{\WW}$ the resulting matrix obtained with the greedy procedure. The step-by-step process is shown in the blue box (Algorithm 1). Besides optimizing the conditioning of the matrix in the normal equations, this approach for the selection of additional modes has the physical meaning of adding experimental measurements at locations where the most relevant dynamics can be observed. 

To summarize the overall cost of our procedure, we have gathered the dominant computational costs for every step within the pipeline, as shown in Table \ref{tab:cost}. The heavy computational burden of the first three steps is kept at the training (or offline) stage, so that whenever a reconstruction must be carried out, only a small number of computations are left, scaling with the number of modes of the SVD. All the steps have behaved well for large scale problems, with the exception of the SVD construction, which could lead to excessive RAM memory consumption. In such cases, the classical method is replaced by the randomized SVD procedure \cite{romor2025dataassimilationperformedrobust,rSVD2024}, where a stochastic sampling of the column space for the snapshot matrix is carried out, leading to fastest computations, in case it is needed.

\noindent \textbf{Remark:} The greedy approach inherits the mathematical analysis of the full method. Let us denote
$$
\beta(m) = \inf_{c \in \bR^n} \left(\frac{\GG c}{\ROM c} \right),
$$
which controls the stability of the recovery procedure for a given number $m$ of measurements. Observe that our greedy algorithm, by construction (see step 2 of Algorithm 1), is designed so that, if we pick $m^*$ measurements with our method, then
$$
\beta(m^*) \geq \beta(m_{\text{sub}}),
$$
for any other selection of sub-optimal measurements $m_{\text{sub}}$, thus leading to the best possible and robust reconstruction for a given number of available measurements. In other words, selecting the sensor that gives the largest possible $\beta(m)$ implies selecting the sensor that gives the smallest possible error according to bound \eqref{eq:apriori}.

\begin{table}[!htbp]
\caption{Computational cost of the full pipeline. The free constant $\alpha$ could be up to 2 depending on the iterative scheme used to invert the system of equations and the preconditioning quality. Only the reconstruction and assembly step is done whenever new measurements are received. Taking into account that $n$ is significantly smaller than the full order model degrees of freedom $\cN$, we end up with near real-time reconstructions.}
\label{tab:cost}
\begin{tabular}{@{}cccccc@{}}
\toprule
                                                                          & \begin{tabular}[c]{@{}c@{}}Solution \\ sampling \\ via FVM\end{tabular} & \begin{tabular}[c]{@{}c@{}}Training\\  SVD\end{tabular} & \begin{tabular}[c]{@{}c@{}}Computation of \\ observation \\ matrix $\textbf{W}$\end{tabular} & \begin{tabular}[c]{@{}c@{}}Reconstruction with\\  $n \ll \mathcal{N}$ modes\end{tabular} & \begin{tabular}[c]{@{}c@{}}Assembly \\ of $u^*$\end{tabular} \\ \midrule
\begin{tabular}[c]{@{}c@{}}Dominant \\ computational \\ cost\end{tabular} & $\mathcal{O}(k \mathcal{N}^\alpha)$                                     & $\mathcal{O}(k^2 \mathcal{N})$                          & $\mathcal{O}(m \mathcal{N}^\alpha)$                                                          & $\mathcal{O}(n^\alpha)$                                                                  & $\mathcal{O}(n \mathcal{N})$                                 \\ \bottomrule
\end{tabular}
\end{table}

\section{Numerical examples}\label{Sec_NumEX}

This section focuses on evaluating the proposed inverse methodology in various settings.  First, we validate the reconstruction of the full temperature field from synthetic thermo-camera measurements acquired over the whole computational domain. Secondly, we show how the method can be used to extrapolate the temperature field in the salmon slice using only synthetic measurements within the airflow domain. Thirdly, we apply the proposed Greedy algorithm to determine optimal sensor placement, and verify its robustness as the number of measurements diminishes.

\begin{figure}
    \centering
    \includegraphics[width=1\linewidth]{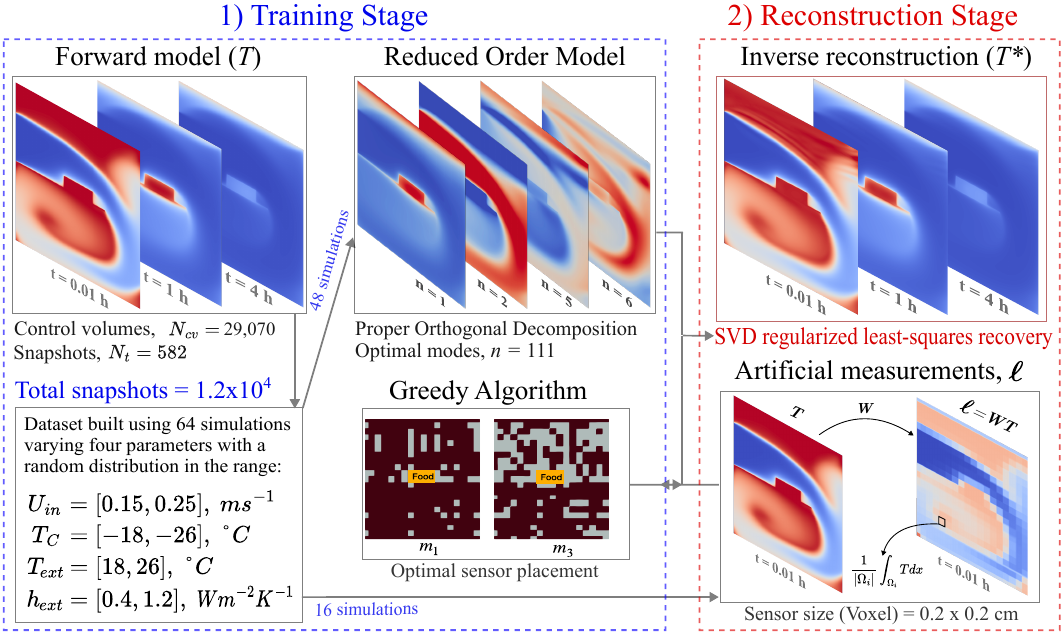}
    \caption{Flow diagram of the proposed inverse framework used to reconstruct the temperatures. A training stage is performed to both, sample the solution space with several numerical simulations of the system dynamics \eqref{eq:govfood}, and to extract relevant flow features with the ROM basis $\ROM$. Next, temperature measurements are received and assimilated to find their best fit with respect to the column space of $\ROM$}
    \label{fig:pipeline}
\end{figure}

The pipeline of the reduced-order inverse estimation framework proposed in this work is illustrated in Fig.~\ref{fig:pipeline}. In the training stage,  we construct a five-dimensional parameter set by gathering snapshots from 64 simulations of the direct model. The resulting time-dependent dataset is obtained by randomly varying the boundary conditions, i.e., $U_{in}$, $T_{C}$, $T_{ext}$, and $h_{ext}$, within a fixed range by using a uniform random distribution, as presented in the lower-left corner of Fig.~\ref{fig:pipeline}. The magnitude of these variation ranges is consistent with real operating conditions of convective freezers and was estimated from the literature review and technical information from previous works \cite{RIVERA2024110558,takhar2017-review1}. Then, we train the reduced order space by using an SVD with 48 simulations from the previous dataset, resulting in a total of 1.2$\times\dot10^{4}$ snapshots. The first 10 modes of the resulting bases are depicted in Fig.~\ref{fig:error_general}.a.
Thereafter, we manufacture artificial low-quality measurements of the temperature fields from the direct model by degrading the resolution in sensors or pixels of size $2\times 2$ cm$^2$ (see also Fig.~\ref{fig:pipeline}). Finally, using a Greedy algorithm, we find optimal sensor placements based on \textit{a-priori} error optimization.
In the reconstruction stage,  we recover the time-dependent temperature fields from 16 synthetic measurements using direct simulations different from those within the training dataset, using the ROM-based least-squares approach described in the previous section.

For the computational implementation of the numerical experiments, for the finite volume scheme solution of the forward problem, we have used an in-house library called MVF2D \cite{MVF2D}. For the full ROM and inverse pipeline, we have used the software MAD (\cite{galarceThesis, MAD}, chapter 5), used successfully for forward and inverse problems in previous works such as \cite{galarce2023_MRE,GLM2022}. The hardware used in the computations for this section were done in a middle-range server using a 48 cores AMD Epyc processor, and 128GB of RAM memory.

\subsection{Reconstruction of temperature fields from thermocamera measurements} \label{sec_thermocamera}

\begin{figure}[h]
    \centering
    \includegraphics[width=1\linewidth]{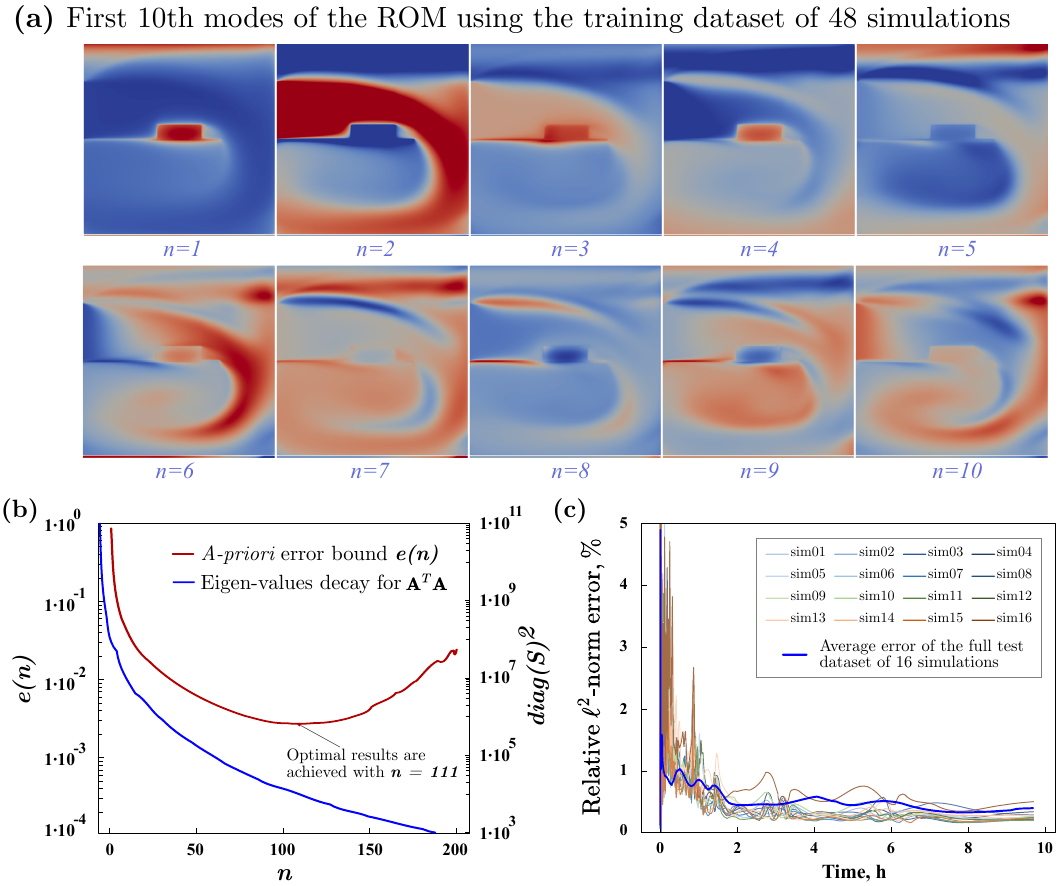}
    \caption{(a) First 10 SVD basis functions. (b) Eigenvalue decay and \textit{a-priori} error as the number of the ROM modes $n$ increases. (c) Evolution of the relative $\ell^2$ norm error for the 16 simulations out of the training set. We observe an agreement between the theoretical \textit{a-priori} bound from \eqref{eq:apriori} and the resulting empirical average error for the full test data-set, reaching a 5 \% peak at the freezing start and staying below the 1 \% until the process ending.}
    \label{fig:error_general}
\end{figure}

The first test considers full-domain synthetic measurements mimicking a thermocamera placed over the two-dimensional mid-plane of the freezer. This setup results in an image composed of $m=361$ pixels or sensors.
Having this setup fixed for the measurement vectors $\ell_i$, one can evaluate the matrix conditioning of the normal equations \eqref{eq:normal_eqs} summarized in \eqref{eq:apriori}. 
The decay of the eigenvalues for the covariance matrix associated with the snapshot matrix is shown in Fig.~\ref{fig:error_general}. The values and decay tendency of the eigenvalues are related to the ROM's ability to approximate the original model's physics. 
The rapid decay of the values suggests, in particular, that the SVD provides a good approximation of the temperature field in the phase change phenomenon, also with a limited number of modes. 

The indicator $e(n)$ from \eqref{eq:apriori} for the fixed measurement setup is shown in Fig.~\ref{fig:error_general}$(a)$. Since this quantity is related to an upper bound for the average error when computing the reconstruction, we infer an optimal dimension for the reduced-order model of $n=111$.
The curve also clearly shows how adding more modes affects the matrix conditioning, lowering the quality of the expected result due to the reduced observability of the ROM when $n$ approaches $m$.

\begin{figure}[h]
    \centering
    \includegraphics[width=1\textwidth]{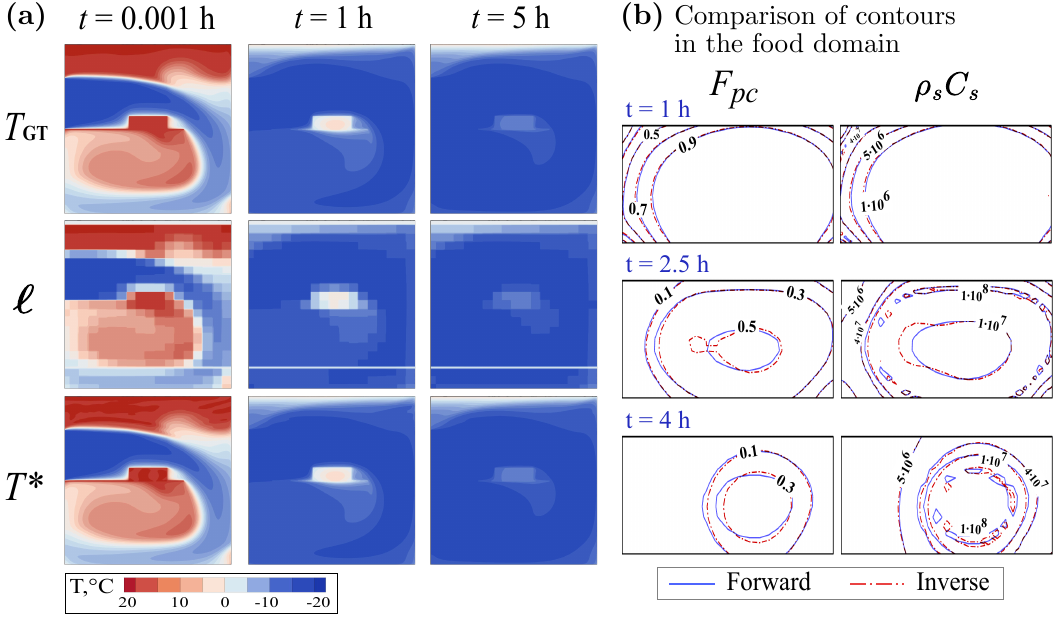} 
    \caption{(a) Temperature field evolution of the ground-truth from the direct model, the full-domain measurements, and the reconstruction procedure at three characteristic snapshots. (b) Comparison of contours of liquid fraction and volumetric heat capacity between the forward and inverse results.}
    \label{fig:rec_minimal}
\end{figure}

A comparison of the temperature fields at three different snapshots is provided in Fig.~\ref{fig:rec_minimal}, showing the ground truth $T_{\text{GT}}$ (not included in the reduced-order model training set), the measurements $\ell = \WW^T T_{\text{GT}}$, and the reconstructions $T^*$. As a benchmark criteria, we introduce the relative $\ell^2$ error:
\begin{equation}
  \text{Relative $\ell^2$-norm error (\%)} = \frac{\norm{T_{\text{GT}} - T^*}}{\norm{T_{\text{GT}}}} \times 100.
  \label{eq:error_apost}
\end{equation}
which we evaluate for every test case. The average result for the 16 simulations is depicted with respect to time in Fig.~\ref{fig:error_general}$(b)$, where we can conclude an overall good reconstruction quality during the whole freezing process, thus validating our framework at a proof-of-concept stage. The reconstruction process and quality depends on how similar the unknown state is to any combination of those in the snapshots dataset. In our case, we train a robust reduced-order model varying the boundary conditions of the forward problem, so that we can ensure a good sensitivity to this parameter.

In addition, we do test exhaustively the reconstruction algorithm using a single, randomly chosen, parameter configuration, with the following values:  $U_{in}=0.2\ ms^{-1}$, $T_{C}=-21.7\ \degree C$, $T_{ext}=21.3\ \degree C$, and $h_{ext}=0.88\ Wm^{-2}K^{-1}$. The results are depicted in Fig.~\ref{fig:rec_minimal}$(a)$ for 3 relevant time-steps of the full field, and in Fig.~\ref{fig:rec_minimal}$(b)$ for a comparison of the liquid fraction contours.
We observe a good qualitative agreement between the ground truth and the reconstructed field during the whole freezing process, especially within the food domain, which is of utmost interest for the targeted application. The field reconstruction indicates that the turbulent airflow recirculation through the cabinet is properly described. 
Moreover, the congruence in the concentration of freezing fronts within the food demonstrates that the phase change process is adequately represented. However, slight differences are noticeable in the upper region of the cabinet, where the airflow is more chaotic due to turbulent fluctuations and the stratification by buoyant forces.

In the food freezing process, the local estimation of heat transfer quantities is of great relevance given that the final quality of frozen food depends mainly on the local freezing rate and the surface heat fluxes it receives. Therefore, in this Section we compare the evolution of two quantities of interest and their local relative error associated only with fixed points within the food for thermocamera measurements. Later, in Section \ref{sec_optimalsensor} we compare the full-domain errors in the form of the l2-norm. The local errors are calculated relative to the time-average values of the quantity of interest $\Phi$ at the point $P$, defined as:
\begin{equation}
    \Phi_{rel,error} (\%) =\frac{|\Phi_{P}-\Phi^{*}_{P}|}{|\overline{\Phi}_{P}|}\times 100.
    \label{eq:localerror}
\end{equation}

\begin{figure}[h]
    \centering
    \includegraphics[width=1\linewidth]{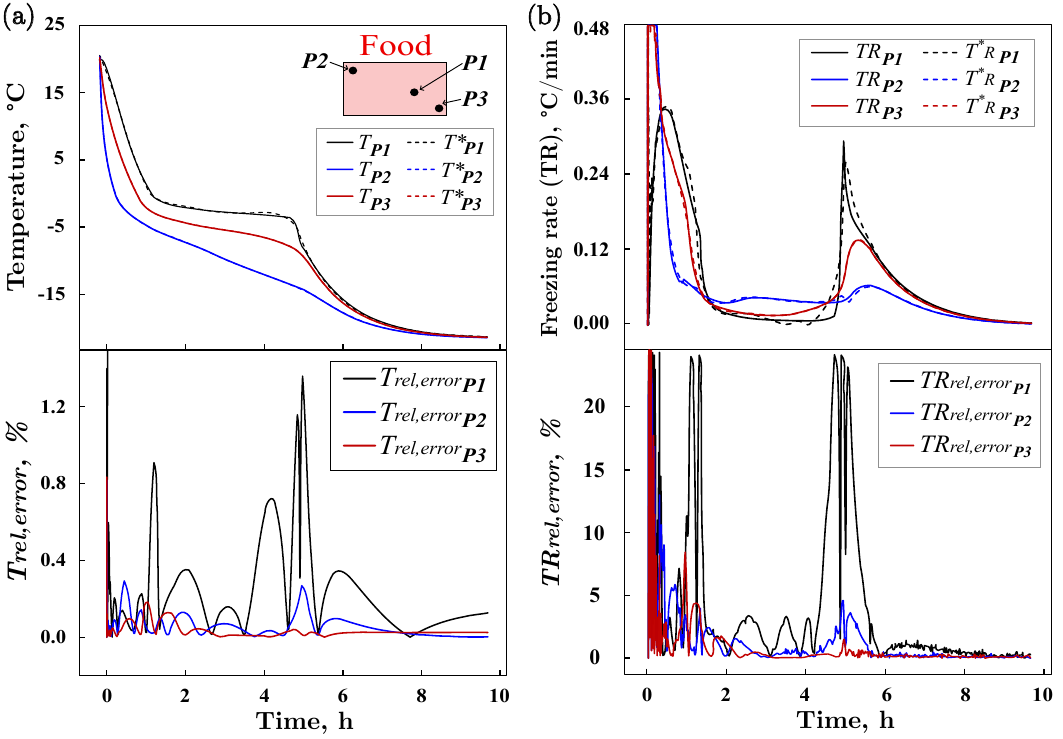}
    \caption{Reconstruction comparison over time for both the temperature (a) and freezing rates (b) at 3 control points. In any of the points, the reconstruction relative error is kept below 2 \% for the temperature.}
    \label{fig:full_domain_exp_errors}
\end{figure}

Fig.~\ref{fig:full_domain_exp_errors} depicts the local evolution of temperature and freezing rate at three fixed points inside the food domain, as well as presents their local reconstruction relative errors. The selected points correspond to the approximate freezing center P1 at $(x,y) = (5.2,\ 1.5)$ cm, and the two opposite corners where freezing rates are expected to be higher at $(x,y) = (0.2,\ 2.8)$ cm for P2, and at $(x,y) = (7.8,\ 0.2)$ cm for P3.
The results indicate that the proposed framework achieves a remarkable precision in reconstructing the evolution of temperatures in the whole food domain since the differences between both models are almost indistinguishable (see Fig.~\ref{fig:full_domain_exp_errors}(a)). However, in the evolution of the freezing rate, the differences between the direct problem and its reconstruction are more noticeable, especially at the freezing center point (P1). The reconstruction method is designed for the field itself, and not its derivatives, so it is still relevant how the errors are kept below a reasonable threshold, as no direct control over this quantity is imposed in the recovery process.

By analyzing the time-dependent relative errors of temperature and freezing rate, we find that the control point in the freezing center presents larger errors than in the food boundaries for the whole process. We must consider that the freezing center represents the point where the freezing fronts concentrate after passing through the entire food domain, so it accumulates all the numerical errors in space and time of the nonlinear food phase-change phenomenon. Therefore, it can be considered that this area is prone to greater numerical errors, which could explain the lower quality reconstruction. Moreover, the local errors increase mainly in two characteristic periods of the freezing process: first, at the very beginning of the initial cooling stage (before 0.1 h), the relative error reaches 1.65\%, and then, at the end of the freezing stage (from 4 to 6 h), where the relative errors get to 1.35\%.
Finally, the relative errors are directly correlated with the freezing rate, which is expected given that the abrupt changes in temperature require an improved temporal discretization. Nevertheless, the relative errors can also be related to the chaotic initial cooling stage of the cabinet which causes slight discrepancies, as illustrated in Fig.~\ref{fig:rec_minimal}, and with the end of the phase change process that accumulates the local error previously discussed.

\subsection{Reconstruction of the validation test with experimental data} \label{sec_recoBeef}

This section discusses a test case to validate the reconstruction capabilities against the experimental data-set of section \ref{sec_validationfvm}, so that we verify that the optimization problem is capable of capturing, with a minimal configuration, the temperature evolution within the control points on the food, where real data is available. To do so, we keep the ROM configuration at its minimum, using the forward simulation from section \ref{sec_validationfvm} to build up the reduced basis $\ROM$, therefore using a single parameter for this task: time. As with the previous example, we consider square pixels with local averages of the temperature. Following the a-priori analysis made for the previous example, we select the optimal ROM dimension to $n=86$, looking at the minimal value of \eqref{eq:apriori}.

\begin{figure}[h]
    \centering
    \includegraphics[width=1\linewidth]{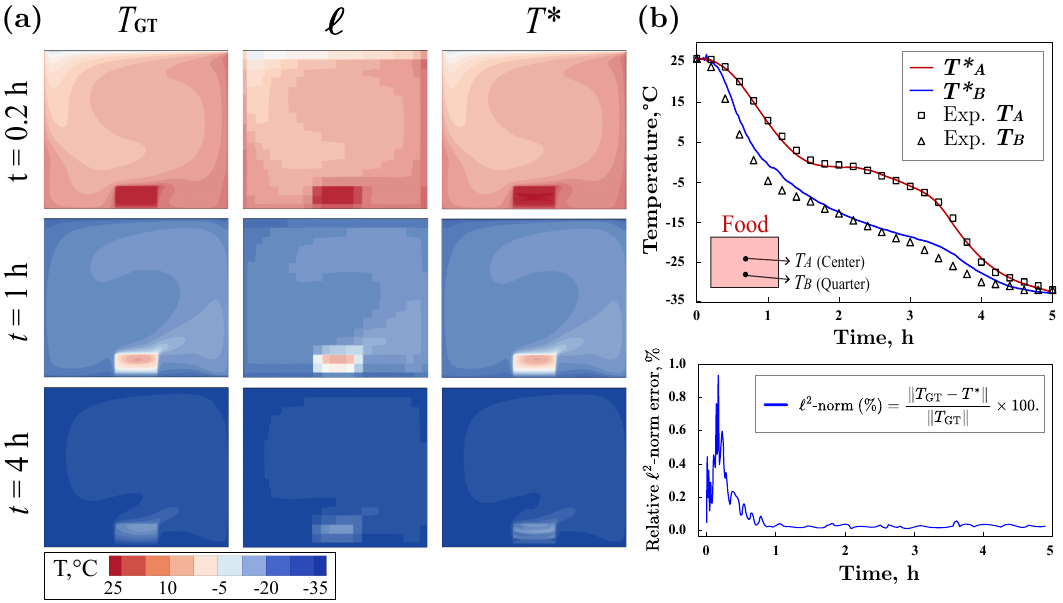}
    \caption{Validation of the inverse framework using test-case with experimental available measurements. (a) The method is able to discover the full temperature field within the working domain from low resolution data, and (b) to fit the empirical data at two control points.}
    \label{fig:validation_inverso}
\end{figure}

The full reconstruction process for this  simple test case is depicted in Figure \ref{fig:validation_inverso}, where we observe an excellent recovery as time goes by, with a peak reaching the 1 $\%$ near the freezing starting stage. In addition, the reconstruction behavior is also in a reasonable range for the three control points with available measurements, leading to a prediction which is as accurate as what the forward model provides.

\subsection{Food temperature reconstruction with extrapolated data from the airflow domain} \label{sec_foodextrapolation}

\begin{figure}[h]
    \centering
    \includegraphics[width=1\linewidth]{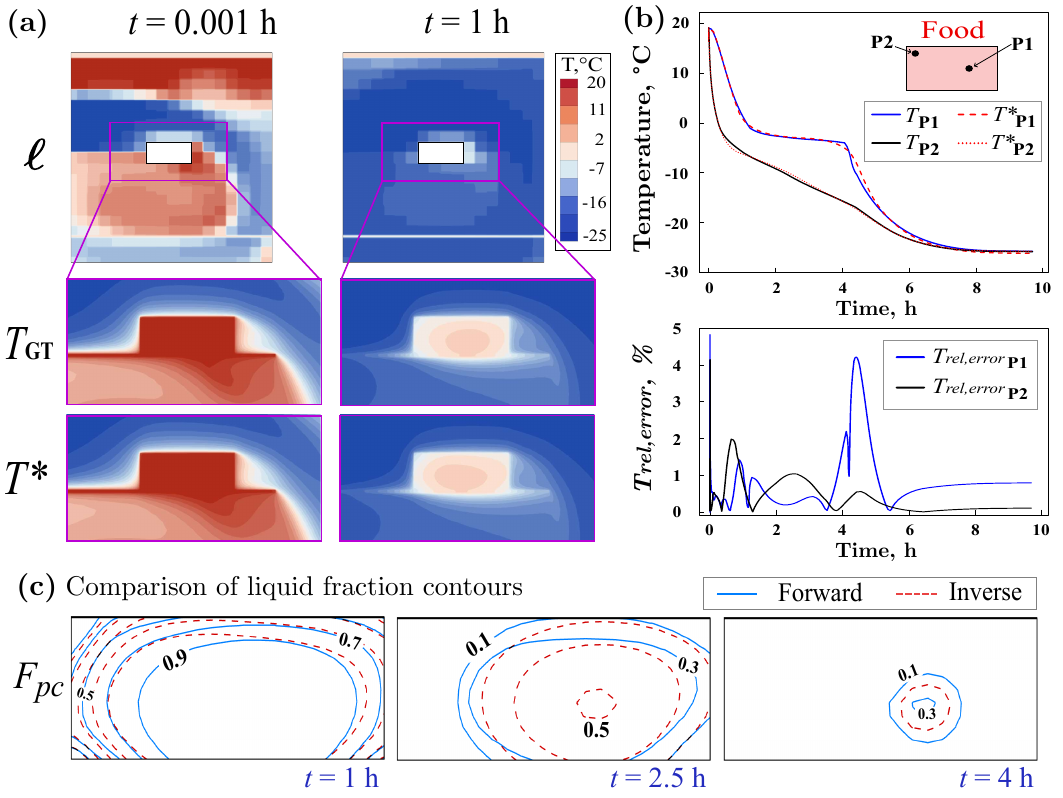}
    \caption{Depiction of temperature recovery using measurements which are blind to the food region, thus allowing us to extrapolate the field beyond the observed region.(a) Fields comparison in the full domain and in the region near the food at two snapshots. (b) Local temperature evolution and its relative error at two fixed points of interest. (c) Comparison of liquid fraction contours between the forward and inverse models at three snapshots.}
    \label{fig:nofood_figs}
\end{figure}

In this section, we consider an extension of the proposed approach to the case when temperature measurements are not available within the food domain $\Omega_{s}$. Formally, we restrict our observation to the airflow region, thus addressing an extrapolation problem for the temperature field in the food. We choose a different set of parameters (from the 16 synthetic dataset) than those used in the previous Section \ref{sec_thermocamera}, with the following values: $U_{in}=0.2$ m$s^{-1}$, $T_{C}=-25.8 \degree$C, $T_{ext}=18.6\ \degree$ C, and $h_{ext}=0.55$ Wm$^{-2}$K$^{-1}$. The rest of the experiment setup is analogous to the previous section, except for the number of sensors, which is set to $m=352$, excluding the pixels in the flood domain. 
Fig.~\ref{fig:nofood_figs}$(a)$ compares the truncated domain of measurements and zoomed sections of the ground-truth and the reconstructed fields for two snapshots, Fig.~\ref{fig:nofood_figs}$(b)$ presents the local evolution and its relative error at two fixed points within the food domain, and  Fig.~\ref{fig:nofood_figs}$(c)$ compares the contours of liquid fraction at three snapshots between the forward and inverse results.

As expected, the quality of the overall reconstruction decreases as the number of sensors is reduced. Despite this, the reconstruction of the temperature fields still achieves great fidelity with respect to the ground-truth; Only slight spatial differences are observed in the freezing fronts within the food.
Similar to the previous section, the quantitative differences in temperature evolution between the ground-truth and the reconstruction are hardly distinguishable, except for a short lapse at the end of the freezing stage. Also, the evolution of relative errors shows a similar tendency with two main periods where the errors peak: at the beginning of the initial cooling stage, and at the end of the freezing stage. Nevertheless, when comparing the magnitude of the relative errors with those of Fig.~\ref{fig:full_domain_exp_errors}, is noticeable that the average errors almost double throughout the entire freezing process, reaching maximum values of 4.8\% at the beginning and 4.2\% at the end of the freezing stage.
The contours of the liquid fraction follow the same pattern described by the isotherms: at the first instant of time, the freezing fronts start to advance from the lateral and upper surfaces, but the latent heat contours do not yet cover large regions of the food; then, at $t=1.5$ h, the freezing fronts concentrate in the center with a slight inclination towards the right lateral surface, while the specific heat contours show clear peaks due to the latent heat; finally, at $t=4$ h, there remains a small area that is not yet frozen, which still concentrates a latent heat peak.
Regarding the reconstruction quality, the results generally indicate that the contours advance at a similar speed but with a slight phase shift throughout the process. In addition, a flatter distortion in the shape of the contours is observed, especially at $t=1.5$ h. Considering the long duration of the process ($t \approx 10$ h) and the nonlinear nature of the temperature and derived quantities evolution, the reconstruction procedure achieves excellent physical representation concerning the ground-truth from the forward model.

This result is very relevant, as it allows us to conclude that the food temperature fields can be reliably estimated without intervening in food, foreseeing this numerical tool as a powerful instrument to complement physical sensors in the food industry.

\subsection{Temperature reconstruction with an optimal sensor placement algorithm} \label{sec_optimalsensor}

\begin{figure}[h]
    \centering
    \includegraphics[width=1\linewidth]{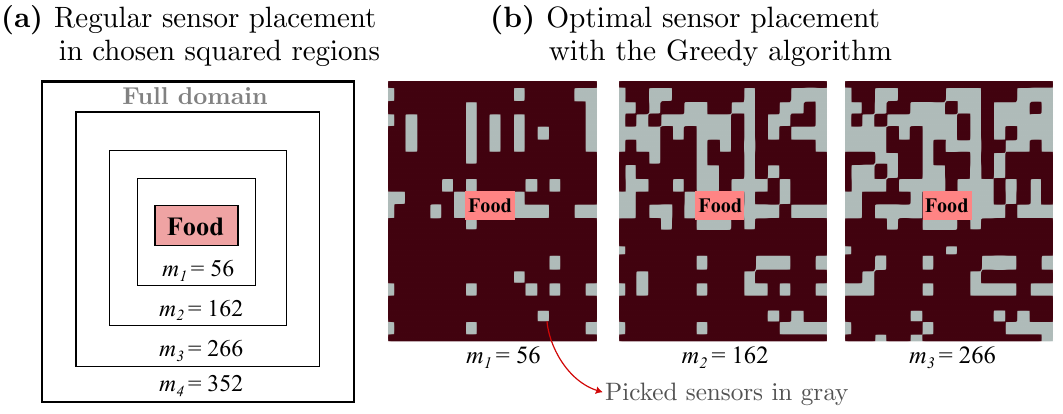}
    \caption{Greedy sensor placement. We compare the temperature reconstruction using a fixed number of measurements ($m_1 = 56$, $m_2=162$, $m_3=266$ and $m_4=352$), always discarding the food chunk. We observe how the methodology provides a sparse spatial sampling of the domain based on the pre-built reduced order model.}
    \label{fig:sensor_placement_congreedy}
\end{figure}

\begin{figure}[h]
    \centering
    \includegraphics[width=1\linewidth]{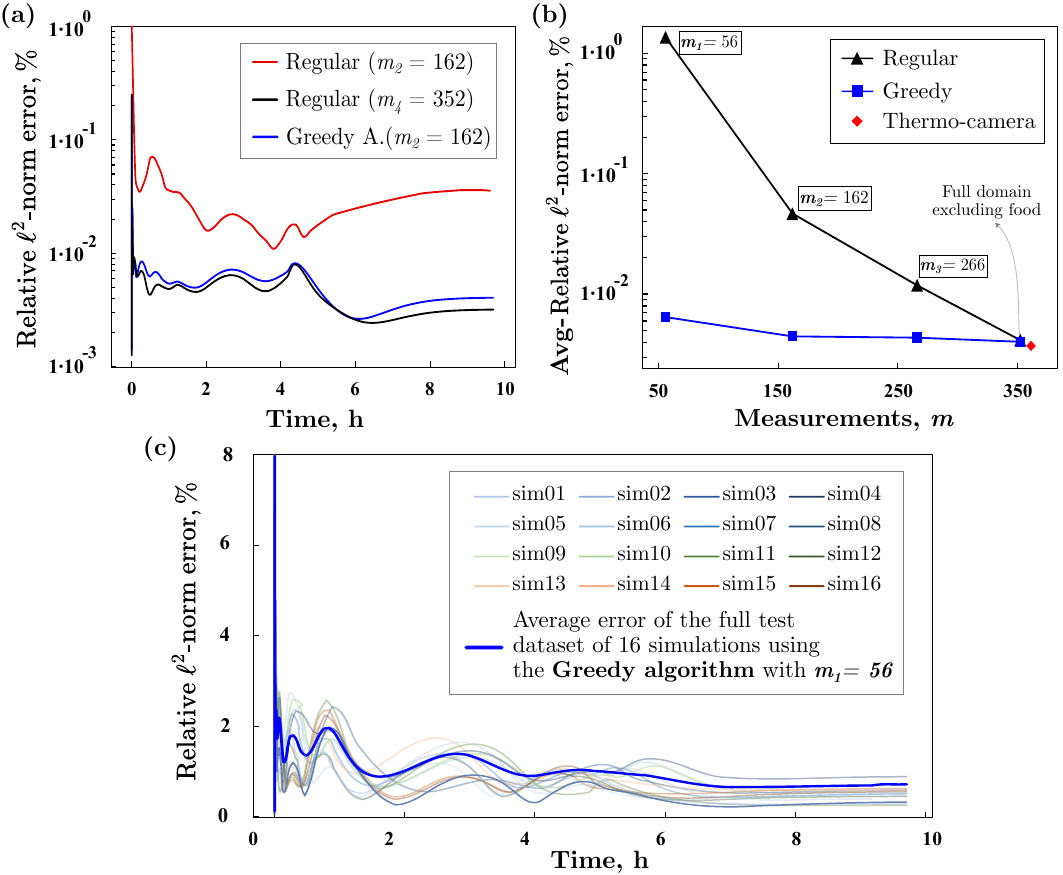}
    \caption{Numerical study with greedy algorithm for sensor placement. We observe how the method selects sensors to ensure better accuracy in (a). The method delivers a good behavior against the number of measurements, showing much better convergence rate in (b) with respect to the regular sensor grid selection. In (c), we see this result systematically tested against 16 different simulations, showing good sensitivity against a number of freezing configurations including different boundary conditions.}
    \label{fig:comp_greedy}
\end{figure}

In this section, we test the optimal sensor placement algorithm proposed in Section \ref{sec:greedy} using the same set of parameters for the boundary conditions used in Section \ref{sec_foodextrapolation}. We compare the reconstruction output of this algorithm with a regularly chosen observation region, as depicted in Fig.~\ref{fig:comp_greedy}. Nested squared observation regions are taken, excluding the food region, leading to 4 sensor set-ups, namely with 56, 162, 266, and 352 (full domain excluding food) measurements. In contrast, a sparser sensor selection is made when using the 3-step algorithms, as shown in Fig.~\ref{fig:sensor_placement_congreedy}. Interestingly, the method tends to privilege both the food surroundings and the upper part of the computational domain, corresponding to the region where the turbulent airflow interacts with the food and with the stratified air in the upper part of the cabinet.

To evaluate and compare the performance of the proposed Greedy algorithm against regular sensor placement, we compute the relative l2-norm error of temperature defined in Eq. \eqref{eq:error_apost}. The time-dependent and time-averaged values of the relative l2-norm error between regular and optimal sensor placements with different amount of measurements are compared in Fig.~\ref{fig:comp_greedy} $(b)$, and \ref{fig:comp_greedy}$(c)$, respectively.

The results demonstrates the robustness of the Greedy algorithms also when decreasing the number of measurements. We observe how properly placing the sensors improves the quality of the reconstructions during the entire freezing process, achieving errors close to 1\% over the whole time-dependent dynamics. Fig. \ref{fig:comp_greedy}$(a)$ shows that the optimal placement of 162 sensors achieves a similar error (in $\ell^2$-norm) over time as the regular placement with almost twice as many sensors ($m$ = 352). 
Moreover, when comparing the time-averaged errors for a given number of sensors (Fig. \ref{fig:comp_greedy}, the Greedy algorithm outperforms the regular placement, especially for a small number of measurements (such as for $m = 52$), reaching similar time-average errors that those obtained with full-domain measurements (thermo-camera). We observe a good behavior against the number of sensors in Fig. \ref{fig:comp_greedy}$(b)$, where the method manages to keep the error low despite the diminished amount of measurements and we also observer a good sensitivity for different freezing conditions in Fig. \ref{fig:comp_greedy}$(c)$, where several boundary conditions and physical parameters are tested out, delivering an average error in time no large than 2 \%.

\begin{figure}[h]
    \centering
    \includegraphics[width=1\linewidth]{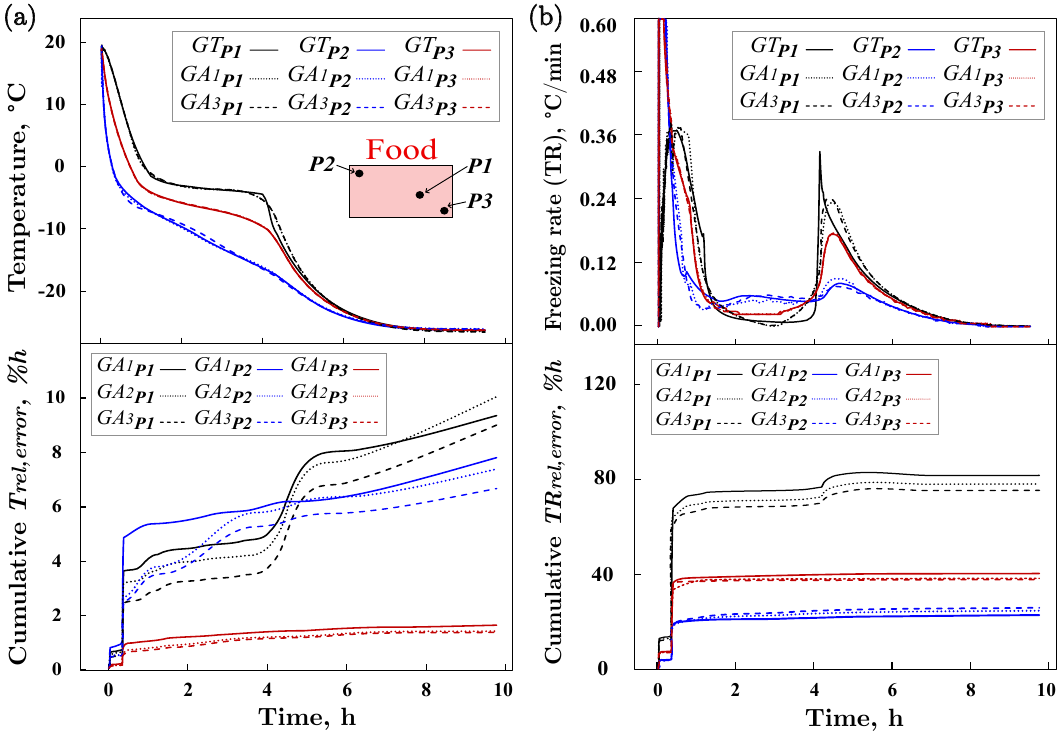}
    \caption{Temperature evolution (a), freezing rates (b), and cumulative errors (c,d), for the food chunk in 3 control points using the greedy reconstruction algorithm. We observe, as expected, an excellent behavior for the temperature reconstruction within the food piece, whereas the freezing rate exhibits an admittedly larger error, as no control is made over the field derivatives during the reconstruction process.}
    \label{fig:fwdvsgreedy_graphs1}
\end{figure}

Similar to the previous sections, we evaluate the reconstruction quality for two different local quantities of interest. This time, we compare the ground-truth with the Greedy algorithm for the three set of sensors, i.e. $GA_{1}$ for $m = 52$, $GA_{2}$ for $m = 162$, and $GA_{3}$ for $m = 266$, for the local evolution of temperature in Fig. \ref{fig:fwdvsgreedy_graphs1}$(a)$, and for the local freezing rate in Fig. \ref{fig:fwdvsgreedy_graphs1}$(b)$. Additionally, we compute the accumulative relative error over time by integrating in time the Equation \ref{eq:localerror}. Additionally, we should note that the sensor location is magnitude-independent, meaning that we only need as an input a set of location candidates, making the greedy search only needed once for all the test cases.

From Fig. \ref{fig:fwdvsgreedy_graphs1}$(a)$ and $(b)$, we can highlight that the local temperature evolution and the freezing rate are well reproduced with the Greedy strategy using $GA_{1}$ and $GA_{3}$ in points P1 and P3. For point P2, $GA_{3}$ depicts a better reconstruction quality, particularly in the freezing rate during the phase-changing stage. Regarding the accumulative relative error over time, we can highlight that half of the total temperature error and two-thirds of the total freezing rate error accumulated at the end of the process comes from the initial chaotic stage, which is associated with the unsteady turbulent fluid dynamics and the rapid cooling of the food boundaries. At this stage, the time derivative of the reconstructed field is of relevance. Despite this, the method reconstructs the temperature field very well and even achieves a good definition of its temporal evolution, which is not trivial since the proposed algorithm does not incorporate any condition regarding the time derivative of the reconstructed field.

\section{Conclusions}\label{Sec_Conclusion}

This study has addressed the inverse estimation of temperature fields in food-freezing processes, a critical factor in ensuring product quality in the food industry. Numerical simulations of thermally coupled flow systems, including liquid-to-solid phase change and turbulence modeling, have been used to generate a synthetic dataset and create a physics-informed reduced-dimension space to regularize the optimization problem for the state estimation. The approach effectively estimates the time-dependent temperature distribution inside the food from sparse sensor data, even when observations exclude the food domain. 

The optimal placement of temperature sensors used for reconstruction was determined using a new greedy algorithm based on available \textit{a priori} estimates, which enhances accuracy and minimizes the number of sensors required. Furthermore, we designed a real-time reconstruction approach, which involves solving a linear system of order $\sim 10^2$, applicable to any new dataset. The novelty of using this approach goes beyond the fact that this strategy had not yet been employed in food freezing applications; most importantly, our approach only requires the ROM and the sensor positions, but not the values of the measurements themselves, enabling the sensor location computation a-priori and independent of the data being assimilated at a given iteration. The Greedy approach has been successfully validated for defining sensor locations that maximized data representativeness, offering a cost-effective solution for complex temperature monitoring applications in food preservation. 

Future work will explore further refinements in the reconstruction algorithm, including model bias \cite{HAIK2023115868}, and/or measurement bias \cite{GMC2024}, control over the time derivative of the reconstructed field, as well as real-time estimation techniques to advance food safety protocols. From the physical viewpoint, future extensions will include Large-Eddy Simulation turbulent models instead of URANS strategies and the temperature estimation of liquid-to-solid scenarios where the natural convection inside the freezing food cannot be despised. The application of the proposed framework for multiphase food freezing or 3D simulations is also a path to strengthen the credibility of the model in the food industry.

\section*{Acknowledgments}

E.~Castillo acknowledges the support provided by ANID-Chile throughout  the project FONDECYT 1210156. D.R.Q.~Pacheco acknowledges funding by the Federal Ministry of Education and Research (BMBF) and the Ministry of Culture and Science of the German State of North Rhine-Westphalia (MKW) under the Excellence Strategy of the Federal Government and the Länder. F. Galarce acknowledges the support provided by DI-VINCI PUCV funding 039.731/2025. Diego.~Rivera acknowledges the support granted by the DICYT-POSTDOC 052416CD funded by VRIDEI-USACH.



\end{document}